\documentstyle[12pt]{article}
\begin{document}

\def\Empty{}
%
% FLoat placement parameters:
%
\def\bottomfraction{.7}
\def\textfraction{0}
\def\floatpagefraction{.7}
%
% Definitions that use @ :
\catcode`\@=11
% Redefine section heading to get smaller letters (this
% is copied from the original in rep12.sty)

%\def\section{\@startsection {section}{1}{\z@}{-3.5ex plus -1ex minus
%    -.2ex}{2.3ex plus .2ex}{\large\bf}}
\def\section{\@startsection {section}{1}{\z@}{-3.5ex plus-1ex minus
    -.2ex}{2.3ex plus.2ex}{\large\bf}}
\def\subsection{\@startsection{subsection}{2}{\z@}{-3.25ex plus-1ex
    minus-.2ex}{1.5ex plus.2ex}{\large\bf}}
\def\subsubsection{\@startsection{subsubsection}{3}{\z@}{-3.25ex plus
 -1ex minus-.2ex}{1.5ex plus.2ex}{\normalsize\bf}}

% Here is \eqalign from plain TEX, resurrected here because latex
% refuses to make it easy to put equation numbers BETWEEN two
% equations...
\def\eqalign#1{\,\vcenter{\openup\jot\m@th
  \ialign{\strut\hfil$\displaystyle{##}$&$\displaystyle{{}##}$\hfil
        \crcr#1\crcr}}\,}

% Modified description environment with paragraph indentation
% (lifted from report.sty)
\def\mydesc{\list{}{\labelwidth\z@ \itemindent-\leftmargin
\listparindent 1.5em
\let\makelabel\descriptionlabel}}
\let\endmydesc\endlist

% Our own caption style --

\def\fnum@figure{{\small Figure \thefigure}}
\def\fakefigure{\def\@captype{figure}}

\long\def\@makecaption#1#2{
    \vskip 10pt
    \def\FCap{#2} \def\NoCap{\ignorespaces}
    \ifx \FCap\NoCap
       \setbox\@tempboxa\hbox{#1}  % This is to avoid the damn colon.
      \else
       \setbox\@tempboxa\hbox{#1: \small \it #2}
    \fi
    \ifdim \wd\@tempboxa >\hsize   % IF longer than one line:
        \unhbox\@tempboxa\par      %   THEN set as ordinary paragraph.
      \else                        %   ELSE  center.
        \hbox to\hsize{\hfil\box\@tempboxa\hfil}
    \fi}

% Our page heading style:
\pagestyle{headings}
\oddsidemargin 0.5in
\evensidemargin 0.5in
\def\@oddhead{\hbox{}\rightmark \hfil \rm\thepage}% Right heading.
\def\sectionmark#1{\markright {\sc{\ifnum \c@secnumdepth >\z@
      \S\thesection.\hskip 1em\relax \fi #1}}}

\catcode`\@=12

\def\oplabel#1{
  \def\OpArg{#1} \ifx \OpArg\Empty {} \else
        \label{#1}
  \fi}

% Here we define the numbering conventions for theorems, etc.
% They are all numbered together, consecutively in each section.
%
\newtheorem{theoremSt}{Theorem}[section]
\newtheorem{corollarySt}[theoremSt]{Corollary}
\newtheorem{propositionSt}[theoremSt]{Proposition}
\newtheorem{lemmaSt}[theoremSt]{Lemma}
\newtheorem{conjectureSt}[theoremSt]{Conjecture}
\def\NextItem{\refstepcounter{theoremSt} \item[\thetheoremSt]}
% This macro allows us to define a bunch of similar theorem-environments
% with automatic labelling and extra features.
%
\def\MakeStEnv#1{
  \newenvironment{#1}[2]{
  \begin{#1St} \oplabel{##1}%
  \global\def\CrntSt{\thetheoremSt}%
  {\rm ##2}%
}{
  \end{#1St} }
}
\MakeStEnv{theorem}
\MakeStEnv{corollary}
\MakeStEnv{proposition}
\MakeStEnv{lemma}
\MakeStEnv{conjecture}

\def\lref#1{\ref{#1} (#1)}% Prints label as well as number
\newenvironment{proof}[1]{
  \def\PfArg{#1}
  \ifx\PfArg\Empty
        \edef\PfArg{\CrntSt}  \fi
 \startproof{\PfArg}%
}{
  \finishproof{\PfArg}
}
\newcommand{\startproof}[1]{
  \medbreak\mbox{}
  {\it Proof of #1:}%
}
\newcommand{\finishproof}[1]{%
\def\FPArg{#1}\ifx\FPArg\Empty\def\FPArg{\CrntSt}\fi%
\smallbreak\noindent\makebox[\textwidth]{\hfill\fbox{\FPArg}}%
\medbreak\noindent}

\newcommand{\marginwrite}[1]{}
% file macros.tex
%
%%%%%%%%%%%%%%%%%%%%%%%%%%%%%%%%%%%%%%%%%%%%%%%%%%%%%%%%%%%%%%%%%%%%%%%%%%%%%%

%                     ABBREVIATIONS

\def\del{\partial}           % boundary symbol
\def\ol#1{\overline{#1}}     % overline

%%%%%%%%%%%%%%%%%%%%%%%%%%%%%%%%%%%%%%%%%%%%%%%%%%%%%%%%%%%%%%%%%%%%%%%%%%%%%%

%                     MATHEMATICS SYMBOL MACROS

\def\a{\alpha}
\def\b{\beta}
\def\g{\gamma}
\def\d{\delta}

\def\abs#1{\vert#1\vert}     % absolute value
\def\eq{\!=\!}               % improved spacing for equals sign
\def\extendor{\raisebox{0ex}[0in][0in]{\rule{0.5pt}{4ex}}}
\def\mat#1#2#3#4{\pmatrix{#1&#2\cr #3&#4\cr}}  %two-by-two matrix
\def\norm#1{\|#1\|}    % norm  || #1 ||
\def\set#1{\{ #1 \}}         % set braces  {  }
\def\vbar{\hbox{\hbox{$\;\vert\;$}}}% vertical bar with wider
                                    %   spacing, e. g.  {x | x = 2}
%%%%%%%%%%%%%%%%%%%%%%%%%%%%%%%%%%%%%%%%%%%%%%%%%%%%%%%%%%%%%%%%%%%%%%%%%%%%%%

%      MATH WORDS AND SYMBOLS NEEDING ROMAN FONT OR SPECIAL SPACING

\def\Aut{\hbox{\it Aut\/}}
\def\Hom{\hbox{\rm Hom\/}}
\def\Imb{\hbox{\it Imb\/}}
\def\imb{\hbox{\it imb\/}}
\def\Inn{\hbox{\it Inn\/}}
\def\Isom{\hbox{\it Isom\/}}
\def\Out{\hbox{\it Out\/}}

%%%%%%%%%%%%%%%%%%%%%%%%%%%%%%%%%%%%%%%%%%%%%%%%%%%%%%%%%%%%%%%%%%%%%%%%%%%%%%

%              \today      (value is today's date)

\def\today{\ifcase\month\or
   January\or February\or March\or April\or May\or June\or
   July\or August\or September\or October\or November\or December\fi
   \space\number\day, \number\year}

%%%%%%%%%%%%%%%%%%%%%%%%%%%%%%%%%%%%%%%%%%%%%%%%%%%%%%%%%%%%%%%%%%%%%%%%%%%%%%

\newcommand{\degree}{^\circ}
\newcommand{\degrees}{^\circ}
\newcommand{\complexes}{{\bf C}}
\newcommand{\gozto}{\mapsto}
\newcommand{\hyperbolic}{{\bf H}}
\newcommand{\hy}{{\bf H}}
\newcommand{\id}{{\rm id}}
\newcommand{\integers}{{\bf Z}}
\newcommand{\supperp}{^\perp}
\newcommand{\reals}{{\bf R}}
\def\zeemod #1{{\bf Z}_{#1}}
\newcommand{\euclidean}{{\bf E}}
\newcommand{\Euclidean}{{\bf E}}
\def\omicron{\circ}

\def\C{\complexes}
\def\Z{\integers}
\def\R{\reals}
\def\P{{\bf P}}
\def\E{\Euclidean}
\def\H{\hyperbolic}
\def\Cstar{{\bf C_\ast}}
\def\Rstar{{\bf R_\ast}}
\def\SD{{\cal D}}
\def\empset{\phi}

\def\Fr{{\rm Fr\/}}
\def\PGL{{\rm PSL}}
\def\PSL#1{{\rm PSL}{#1}}
\def\PSLtil#1{\til{\rm PSL}{#1}}
\def\SL#1{{\rm SL}{#1}}
\def\GL#1{{\rm GL}{#1}}
\def\CP#1{{\bf CP}^{#1}}
\def\I{{\rm I}}
\def\BDiff{{\it BDiff}}
\def\Diff{{\it Diff}}
\def\diff{{\it diff}}
\def\Isom{{\it Isom}}

\def\union{\cup}
\def\til{\widetilde}
\def\hat{\widehat}
\def\rel{\;\hbox{\rm rel}\;}

\def\uu#1{\hbox{\vtop{\hbox{$#1$}\vskip1.1pt\hrule\vskip1.1pt\hrule}}}

\font\caps=cmcsc10  % all capitals
\def\demo#1{\par\noindent {\caps #1:}}

\def\mapdown#1{\big\downarrow % makes a downarrow with a mathematics
            \rlap            % mode symbol to the right
            {\smash{$\vcenter       % (for use in the "diagram" template)
            {\hbox{$         %
            \scriptstyle#1   %
            $}}$}}}           %

\def\mapright#1{\smash{      % makes a long rightarrow with a mathematics
            \mathop          % mode symbol above it
            {\longrightarrow % (for use in the "diagram" template)
            }\limits^{#1}}}  %

\def\bscr{{\cal B}}
\def\escr{{\cal E}}
\def\gscr{{\cal G}}
\def\hscr{{\cal H}}
\def\kscr{{\cal K}}
\def\b{\uu{b}}
\def\bone{\uu{b_1}}
\def\btwo{\uu{b_2}}
\def\f{\uu{f}}
\def\fone{\uu{f_1}}
\def\ftwo{\uu{f_2}}
\def\s{\uu{\sigma}}
\def\sone{\uu{\sigma_1}}
\def\stwo{\uu{\sigma_2}}
\def\m{\uu{m}}
\def\mone{\uu{m_1}}
\def\mtwo{\uu{m_2}}

\title{Ubiquity of geometric finiteness in mapping
class groups of Haken 3-manifolds}
\author{Sungbok Hong and Darryl McCullough}
\date{{\footnotesize Departent of Mathematics, Korea University,
Seoul 136-701, Korea}
\\
{\footnotesize Department of Mathematics, University of Oklahoma,
Norman, OK 73019, USA}
\\
\bigskip
{\footnotesize July 1, 1997}}
\maketitle

\section{Introduction}
\label{intro}

The mapping class group $\hscr(F)$ of a 2-manifold is the group of
isotopy classes of diffeomorphisms, $\pi_0(\Diff(F))$. As is
well-known, it acts properly discontinuously on a Teichm\"uller space,
which is topologically a Euclidean space (traditionally only the
action of the orientation-preserving classes was considered). This
classical setup was refined by Harer \cite{Harer1,Harer2}, who
found an ideal triangulation of Teichm\"uller space for which
$\hscr(F)$ acts simplicially, then constructed a contractible
simplicial complex in the first barycentric subdivision of this
triangulation which is invariant and has finite quotient under the
action. Consequently, if $\Gamma$ is any torsionfree subgroup of
finite index in $\hscr(F)$ (and such subgroups always exist), the
quotient of the action of $\Gamma$ on this complex is a finite
$K(\Gamma,1)$-complex. That is, $\Gamma$ is a {\it geometrically
finite} group.

For a compact 3-manifold $M$, a great deal is known about $\hscr(M)$
(see for example~\cite{K} or the surveys~\cite{McCPoland,McCKorea}).
In the present paper, we extend the aforementioned property of
torsionfree subgroups of $\hscr(F)$ to the case of Haken 3-manifolds.
To postpone the introduction of boundary patterns, we state here only
a corollary of our main theorem. In the corollary, $\hscr(M \rel W)$
means $\pi_0(\Diff(M\rel W))$, where $\Diff(M\rel W)$ is the group of
diffeomorphisms of $M$ which restrict to the identity on~$W$.

\bigskip
\noindent{\bf Corollary \ref{rel W corollary}} {\em Let $M$ be a Haken
3-manifold and let $W$ be a compact 2-dimensional submanifold of
$\partial M$ such that $\partial M-W$ is incompressible. Then any
torsionfree subgroup of finite index in $\hscr(M \rel W)$ is
geometrically finite.}
\smallskip

\noindent
In this result, $W$, or $\partial M-W$, or both may be empty. For
Haken 3-manifolds, it was already known that $\hscr(M\rel W)$ {\it
contains} a finite-index subgroup which is geometrically finite
\cite{jdg}, so the new information from the Main Theorem is that {\it
every} torsionfree subgroup of finite index in $\hscr(M\rel W)$ is
geometrically finite. It is a longstanding and apparently difficult
open question whether a torsionfree finite extension of a
geometrically finite group must be geometrically finite. Of course, if
this were known, the result from \cite{jdg} would imply our
strengthened version.

Our method of proof is to construct a topological action of
$\hscr(M\rel W)$ on a contractible simplicial complex (whose quotient
is compact), and as we show in section~\ref{agf} this is sufficient to
deduce the geometric finiteness of torsionfree finite-index subgroups.
It would be interesting to give a more direct construction of a
contractible complex, along the lines of those developed by Harer,
admitting a {\it simplicial} action of $\hscr(M\rel W)$ with finite
quotient.

The Kontsevich Conjecture (problem 3.48 in the new version of
R.~Kirby's problem list~\cite{K}) asserts that the classifying space
$\BDiff(M\rel\partial M)$ has the homotopy type of a finite complex
when $M$ is a compact 3-manifold with nonempty boundary. This
conjecture was recently proven for the case of Haken 3-manifolds in
\cite{Hatcher-McCullough}. As observed there, the Kontsevich
Conjecture for a Haken 3-manifold $M$ is equivalent to the assertion
that $\hscr(M\rel\partial M)$ is geometrically finite. In
section~\ref{kconj} we use our results on geometric finiteness to
deduce a generalization of the result from \cite{Hatcher-McCullough}:

\bigskip
\noindent{\bf Theorem \ref{Kontsevich Conjecture classifying spaces}}
{\em Let $M$ be a Haken 3-manifold with incompressible boundary, and
let $F$ be a nonempty compact 2-manifold in $\partial M$ such that
$\partial M-F$ is incompressible. Then $\BDiff(M\rel F)$ has the
homotopy type of a finite complex.}
\smallskip

\noindent The proof makes use of the extension of Nielsen's theorem to
3-manifolds made by Heil and Tollefson~\cite{H-T}.

We will work in the context of 3-manifolds with boundary patterns.
This lends greater generality to the results, and allows us to make
direct use of Johannson's powerful characteristic submanifold theory
for Haken manifolds. As is well-known, the characteristic submanifold
was discovered and exploited independently by Johannson~\cite{Joh} and
Jaco and Shalen~\cite{JS}. We use Johannson's formulation because it
is ideally suited to working with homotopy equivalences and
homeomorphisms of 3-manifolds. In section~\ref{char submfd} we provide
a brief exposition of the portion of Johannson's theory that we will
use. In section~\ref{agf} we introduce a generalization of geometric
finiteness, called (for lack of imagination) {\it almost geometric
finiteness}. Any torsionfree finite-index subgroup of an almost
geometrically finite group is geometrically finite. Using a theorem of
Kamishima, Lee, and Raymond~\cite{KLR}, we prove in
proposition~\ref{algebraic lemma} a key fact: an extension of a
virtually finitely generated abelian group by an almost geometrically
finite group is almost geometrically finite. The rest of the proof
then follows the general approach of~\cite{jdg} to show that
$\hscr(M)$ is an extension of this form, and hence is almost
geometrically finite. It is necessary to work with a bit more
precision than was needed in~\cite{jdg}, since one can no longer evade
difficulties by passing to subgroups of finite index in~$\hscr(M)$. We
also give a new proof that torsionfree subgroups of finite index in
$\hscr(M)$ exist. This was proven in~\cite{jdg} by a very complicated
argument; the new proof uses an algebraic fact from~\cite{M-M} to give
a much shorter and more transparent proof. The final section contains
the application to the Kontsevich Conjecture.

The authors thank BSRI-1422 for support of their collaborative
work, and acknowledge helpful discussions with Kyung-bai Lee and
Leonard Rubin.

\section [The Characteristic Submanifold] {Johannson's characteristic
submanifold theory}
\label{char submfd}\marginwrite{char submfd}

\medskip

We give here a brief review of the basic definitions of Johannson's
formulation of the characteristic submanifold. We refer the reader
to~\cite{Joh} for the original presentation, and also to chapter~2
of~\cite{C-M} for a more extensive expository treatment with a number
of examples.

A {\it boundary pattern} $\uu{m}$ for an $n$-manifold $M$ is a finite
set of compact, connected $(n-1)$-manifolds in $\partial M$, such that
the intersection of any $i$ of them is empty or consists of
$(n-i)$-manifolds. Thus when $n\eq 3$, two elements of the boundary
pattern intersect in a collection of arcs and circles, while if three
elements meet, their intersection consists of a finite collection of
points at which three intersection arcs meet.  It is important in
arguments to distinguish between elements of $\uu{m}$ and the points
of $M$ which lie in them, and we will always be precise in this
distinction. The symbol $\abs{\uu{m}}$ will mean the set of points of
$\partial M$ that lie in some element of $\uu{m}$.  When
$\abs{\uu{m}}\eq \del M$, $\uu{m}$ is said to be {\it complete}.
Provided that $\partial M$ is compact, we define the {\it completion}
of $\uu{m}$ to be the complete boundary pattern $\ol{\uu{m}}$ which is
the union of $\uu{m}$ and the set of components of the closure of
$\partial M-\abs{\uu{m}}$. In particular, the set of boundary
components of $M$ is the completion of the empty boundary
pattern~$\ol{\uu{\emptyset}}$.

Suppose $(X,\uu{x})$ is an admissibly imbedded codimension-zero
submanifold of $(M,\uu{m})$, which is admissibly imbedded in
$(M,\ol{\uu{m}})$. The latter assumption guarantees that $X\cap
\partial M\eq \abs{\uu{x}}$, that $\abs{\uu{x}}$ is contained in the
topological interior of $\abs{\uu{m}}$ in $\partial M$, and that an
element of $\uu{x}$ which does not meet any other element of $\uu{x}$
must be imbedded in the manifold interior of an element of $\uu{m}$.
Let $\uu{x}''$ denote the collection of components of the frontier of
$X$ in $M$. To {\it split $M$ along $X$} means to construct the
manifold with boundary pattern
$(\ol{M-X},\widetilde{\uu{m}}\cup\uu{x}'')$, where the elements of
$\widetilde{\uu{m}}$ are the closures of the components of $F-(X\cap
F)$ for~$F\in\uu{m}$. The boundary pattern
$\widetilde{\uu{m}}\cup\uu{x}''$ is called the {\it proper boundary
pattern} on $\ol{M-X}$.

Maps which respect boundary pattern structures are called admissible.
Precisely, a map $f$ from $(M,\uu m)$ to $(N,\uu n)$ is called {\it
admissible} when $\uu{m}$ is the disjoint union
$$\uu{m}=\coprod_{G\in\uu{n}}\{\hbox{{\rm\ components of\ }}f^{-1}(G)\,\}.$$

\noindent
An admissible map $h\colon (K,\uu{k})\to (X,\uu{x})$, where $K$ is an
arc or a circle and $(X,\uu{x})$ is a $2$-~or $3$-manifold, is called
{\it inessential} if it is admissibly homotopic to a constant map (the
constant map might not be admissible, but all the other maps in the
homotopy must be admissible), otherwise it is called {\it essential.}
A map $f\colon (X,\uu{x})\to (Y,\uu{y})$ between $2$-~or $3$-manifolds
(not necessarily of the same dimension) is called {\it essential} if
for any essential path or loop $h\colon (K,\uu{k})\to (X,\uu{x})$, the
composition $fh\colon (K,\uu{k})\to (Y,\uu{y})$ is essential. Notice
that when $\uu{x}$ is empty, this simply says that $f$ is injective on
fundamental groups.

The group of admissible isotopy classes of admissible diffeomorphisms
from $(M,\uu{m})$ to $(M,\uu{m})$ is denoted by ${\cal H}(M,\uu{m})$.
The classes that preserve the orientation of each component are
indicated by a plus subscript, as in ${\cal H}_+(M,\uu{m})$. The
classes relative to the subset $\abs{\uu{m_1}}$, where
$\mone\subseteq\m$, are denoted by ${\cal
H}(M,\uu{m}\rel\abs{\mone})$.  Suppose $\langle h\rangle\in {\cal
H}(M,\uu{m})$. Since $h^{-1}(\abs{\uu{m}})\eq \abs{\uu{m}}$, $h$ must
carry each free side of $(M,\uu{m})$ diffeomorphically to a free side
of $(M,\uu{m})$. Therefore $h$ is admissible for $(M,\ol{\uu{m}})$.
Thus when working with mapping class groups of manifolds with boundary
patterns, the requirement that the boundary pattern be complete is not
at all restrictive.

An $i$-{\it faced disc} is a 2-disc whose boundary pattern is complete
and has $i$ components. Observe that each element of $\uu{m}$ is
incompressible if and only if whenever $D$ is an admissibly imbedded
$1$-faced disc in $(M,\uu{m})$, the boundary of $D$ bounds a disc in
$\abs{\uu{m}}$ which is contained in a single element of $\uu{m}$. For
most of Johannson's theory, a somewhat stronger condition is needed.
The boundary pattern $\uu{m}$ of a 3-manifold $M$ is called {\it
useful} when the boundary of every admissibly imbedded $i$-faced disc
in $(M,\uu{m})$ with $i\leq3$ bounds a disc $D$ in $\partial M$ such
that $D\cap\,(\cup_{F\in\,\uu{m}}\partial F)$ is the cone on $\partial
D\cap\,(\cup_{F\in\uu{m}}\partial F)$. Notice that
$\ol{\uu{\emptyset}}$ is a useful boundary pattern on $M$ if and only
$\partial M$ is incompressible.

A Seifert fibering on a $3$-manifold $(V,\uu{v})$ with boundary
pattern is called an {\it admissible Seifert fibering} when the
elements of $\uu{v}$ are the preimages of the components of a boundary
pattern of the orbit surface.  Consequently the elements of $\uu{v}$
must be tori or fibered annuli.

Assume that $V$ carries a fixed structure as an $I$-bundle
over~$B$. Each component of the associated $\partial I$-bundle is a
$2$-manifold in $\partial V$, called a {\it lid} of the
$I$-bundle. There are two lids when the bundle is a product, and one
when it is twisted. Let $\uu{b}$ be a boundary pattern on $B$. The
preimages of the elements of $\uu{b}$ form a collection of squares and
annuli in $\partial V$, called the {\it sides} of the $I$-bundle. The
lid or lids, together with the sides, if any, form a boundary pattern
$\uu{v}$ on $V$. When $V$ carries this boundary pattern, the fibering
is called an {\it admissible $I$-fibering} of $V$ as an $I$-bundle
over~$(B,\uu{b})$. We emphasize that for an admissible $I$-fibering,
{\it the lids are always elements of the boundary pattern.}

We are now ready to introduce the characteristic submanifold. An {\it
admissible $I$-bundle} or {\it Seifert fibered space} $(X,\uu{x})$ in
a $3$-manifold $(M,\uu{m})$ is an $I$-bundle or Seifert fibered space
imbedded in $M$ such that the inclusion defines admissible maps
$(X,\uu{x})\to (M,\uu{m})$ and $(X,\uu{x})\to (M,\ol{\uu{m}})$. An
admissible $I$-bundle or Seifert fibered space in $M$ is called {\em
essential} when its frontier is an essential surface in $(M,\uu{m})$.
This implies that the inclusion of $(X,\uu{x})$ into $(M,\uu{m})$ is
an essential map.

We suppose that $(M,\uu m)$ is a Haken 3-manifold with useful boundary
pattern.  A disjoint collection $(\Sigma,\uu{\sigma})$ of essential
admissible $I$-bundles and Seifert fibered spaces is a {\em
characteristic submanifold} for $(M,\uu m)$ if
\begin{enumerate}
\item
$(\Sigma,\uu\sigma)$ is full, i.e. the union of $\Sigma$ with any of
the components of $\ol{M-\Sigma}$ is not a disjoint union of
essential admissible $I$-bundles and Seifert fibered spaces,
\item
(Engulfing property)
any essential admissible $I$-bundle or Seifert fi\-ber\-ed space $(X,\uu x)$
in $(M,\uu{m})$ is admissibly isotopic into~$(\Sigma,\uu\sigma)$, and
\item
(Enclosing property) any essential map $f\colon (T,\uu t)\to (M,\uu
m)$ of a square, annulus, or torus into $(M,\uu{m})$ is
admissibly homotopic to a map with image in~$\Sigma$.
\end{enumerate}
One may combine proposition 9.4, corollaries 10.9 and 10.10 and
theorem 12.5 of \cite{Joh}
to see that every Haken 3-manifold with useful boundary pattern has
a characteristic submanifold, and that the characteristic submanifold
is unique up to admissible isotopy.

A Haken $3$-manifold $(M,\uu{m})$ whose completed boundary pattern
$\ol{\uu{m}}$ is useful is called {\it simple} if every component of
the characteristic submanifold of $(M,\ol{\uu{m}})$ is a regular
neighborhood of an element of $\ol{\uu{m}}$. When the boundary
consists of tori, the interior of a simple 3-manifold admits a
hyperbolic structure of finite volume. Then, by Mostow Rigidity and
Waldhausen's theorem, $\hscr(M,\overline{\uu{\emptyset}})$ is
isomorphic to the group of isometries, hence is finite.  One of the
deep applications of Johannson's theory is the following
generalization of this fact (although it was proven before the
dissemination of Thurston's work on hyperbolic structures of
3-manifolds). It appears as proposition~27.1 of~\cite{Joh}.

\begin{theorem}{FMCGT}
{\rm (Finite Mapping Class Group Theorem):}
Let $(M,\uu{m})$ be a simple $3$-manifold with complete and useful
boundary pattern. Then ${\cal H}(M,\uu{m})$ is finite.
\marginwrite{FMCGT}
\end{theorem}

\section [Almost geometric finiteness] {Almost geometric finiteness}
\label{agf}\marginwrite{agf}

We say that a group $G$ is {\it almost geometrically finite} if it
acts smoothly and properly discontinuously on a contractible
simplicial complex $L$, such that $L/G$ is compact. A subgroup of
finite index in an almost geometrically finite group is almost
geometrically finite. Almost geometrically finite groups are finitely
generated (any group acting properly discontinuously on a connected
locally compact space with compact quotient must be finitely
generated). Every torsionfree subgroup of finite index in an almost
geometrically finite group must be geometrically finite:

\begin{lemma}{torsionfree subgroups are geometrically finite}{} Let
$G$ be an almost geometrically finite group, and let $\Gamma$ be a
torsionfree subgroup of finite index in $G$. Then $\Gamma$ is
geometrically finite.
\marginwrite{torsionfree subgroups are geometrically finite}
\end{lemma}

\begin{proof}{} Let $G$ act properly discontinuously on the
contractible complex $L$ with compact quotient. Since $\Gamma$ has
finite index, $L/\Gamma$ is compact, and since $\Gamma$ is
torsionfree, it acts freely on $L$. Therefore $L/\Gamma$ is locally an
ANR and hence an ANR (see for example theorem~5.4.5 of
\cite{VanMill}). By \cite{West}, a finite-dimensional compact ANR is
homotopy equivalent to a finite simplicial complex, and therefore
$\Gamma$ is geometrically finite.
\end{proof}

\noindent We caution that in general, an almost geometrically finite
group need not contain any geometrically finite subgroups of finite
index: Schneebeli \cite{Schneebeli} constructed extensions $1\to
\Z/k\to E\to Q\to 1$ where $Q$ is geometrically finite but $E$
contains no torsionfree subgroups of finite index.

By lemma~\ref{torsionfree subgroups are geometrically finite}, a group
is geometrically finite if and only if it is almost geometrically
finite and torsionfree. Simple examples of almost geometrically finite
groups are finite groups (acting on a point) and finitely generated
abelian groups (project to the quotient by the torsion subgroup and
act on $\R^n$ by translations). By work of Borel and Serre~\cite{B-S},
many arithmetic groups including $\GL(n,\Z)$ are almost geometrically
finite. Finitely generated virtually free groups are another
interesting example; by~\cite{K-P-S} these are exactly the groups that
act simplicially on trees with finite quotient and finite vertex
stabilizers. In our context, mapping class groups of 2-manifolds are
one of the most important examples:

\begin{lemma}{2-manifold mcg}{}
Let $B$ be a 2-manifold, not necessarily connected, of finite type
with compact boundary, and let $\b$ be a boundary pattern for $B$.
Then $\hscr(B,\b)$ is almost geometrically finite.
\marginwrite{2-manifold mcg}
\end{lemma}

\begin{proof}{} Write $\uu{b}\eq \uu{a}\cup\uu{c}$ where the elements
of $\uu{a}$ are arcs and those of $\uu{c}$ are circles. Let $B_0$ be
obtained by removing all boundary circles that do not contain an
element of $\uu{a}$. Let $P'$ consist of one point from the interior
of each element of $\uu{a}$ and let $P''$ consist of the points that
are the intersection of two different elements of $\uu{a}$. It is not
difficult to see that $\hscr(B,\uu{b})$ is isomorphic to a subgroup
of finite index in $\hscr(B_0,P'\cup P'')$. Thus we are reduced to
proving the lemma for $\hscr(B,P)$, where $P$ is a finite subset of
$\partial B$ that meets every boundary component of~$B$.

We will first prove the lemma under the assumption that $B$ is
connected. Suppose first that $\chi(B)\geq 0$. For $S^2$ with 0, 1, or
2 punctures, $\R\P^2$ with 0 or 1 puncture, $D^2$ with 0 or 1
puncture, the M\"obius band, or the Klein bottle, $\hscr(B,P)$ is
finite, for the annulus it is virtually infinite cyclic, and for the
torus it is $\GL(2,\Z)$. From now on we assume that $\chi(B)<0$.

Suppose first that $P$ is empty, so $B$ is either closed or is a
punctured closed surface. If $B$ is closed and orientable, then
$\hscr(B)$ acts properly discontinuously on Teichm\"uller space.
(Classically, from the viewpoint of conformal surfaces, only the
action of the orientation-preserving classes was considered, but using
the hyperbolic viewpoint as in \cite{F-L-P}, $\hscr(B)$ acts as well.)
This action extends to the Harvey bordification \cite{Harvey1,Harvey2}
of Teichm\"uller space, on which it acts smoothly with compact
quotient. Lifting a smooth triangulation gives the necessary
simplicial structure to conclude that $\hscr(B)$ is almost
geometrically finite. If $B$ is orientable and punctured, then
$\hscr(B)$ acts on the contractible complex $Y^0$ constructed by Harer
\cite{Harer1,Harer2} with compact quotient. The analogues of these for
nonorientable surfaces are given in section~2 of~\cite{jdg}. Suppose
now that $P$ is nonempty. Lemma~1.2 of \cite{Hatcher-McCullough} uses
\cite{HatcherHarer} to construct a contractible complex on which
$\hscr(B,P)$ acts properly discontinuously.

We now assume that $B$ is not connected. Let $B_1,\ldots\,$, $B_n$ be
the components of $B$, and let $P_i\eq P\cap B_i$. Assume first that
each $(B_i,P_i)$ is diffeomorphic to $(B_1,P_1)$. From the connected
case, there is a contractible simplicial complex $L$ such that
$\hscr(B_1,P_1)$ acts properly discontinuously on $L$. Let $L_i$ be a
copy of $L$ on which $\hscr(B_i,P_i)$ acts. Fix diffeomorphisms
$f_i\colon (B_1,P_1)\to (B_i,P_i)$ and equivariant simplicial
isomorphisms $\phi_i\colon L_1\to L_i$. By means of these
identifications, if $\langle h\rangle\in\hscr(B,P)$, and $h$ carries
$B_i$ to $B_j$, we may regard $\langle h\rangle$ as carrying $L_i$ to
$L_j$. For let $h_i\colon B_i\to B_j$ be the restriction of $h$ to
$B_i$. Then let $h_i$ send $L_i$ to $L_j$ by
$\phi_j\phi_i^{-1}\circ\langle f_if_j^{-1}h_i\rangle$. Now we define
an action of $\hscr(B,P)$ on $\prod_{i=1}^n L_i$. Given $h$, define a
permutation $\sigma$ by $h(B_i)\eq B_{\sigma(i)}$. Then, put $\langle
h\rangle(x_1,\ldots,x_n)\eq (\langle
h_{\sigma^{-1}(1)}\rangle(x_{\sigma^{-1}(1)}),\ldots, \langle
h_{\sigma^{-1}(n)}\rangle(x_{\sigma^{-1}(n)}))$. This action is
properly discontinuous with compact quotient, since the finite index
subgroup of $\hscr(B,P)$ that preserves each component acts as a
product of properly discontinuous actions each with compact quotient.

For the general case, partition the components of $(B,P)$ into maximal
subsets such that each subset consists of pairwise diffeomorphic
components. Let $F_1,\ldots\,$, $F_m$ be the unions of the components
in the subsets, and let $Q_i\eq P\cap F_i$. By the previous paragraph,
each $\hscr(F_i,Q_i)$ acts properly discontinuously on a contractible
complex $K_i$. Since $\hscr(B,P)$ is the direct product of the
$\hscr(F_i,Q_i)$, it acts properly discontinuously on the product of
the~$K_i$.
\end{proof}

A useful observation is that if there is a surjective homomorphism
$Q_1\to Q_2$ with finite kernel, and $Q_2$ is almost geometrically
finite, then $Q_1$ is also almost geometrically finite. In our proof that
mapping class groups of Haken 3-manifolds are almost geometrically
finite, we will use the following generalization of this observation,
which is a direct consequence of a theorem of Kami\-shi\-ma, Lee, and
Raymond~\cite{KLR}:

\begin{proposition}{algebraic lemma}{} Let $1\to V\to G\to Q\to 1$ be
an exact sequence of groups, such that $V$ contains a finitely
generated abelian group of finite index and $Q$ is almost
geometrically finite. Then $G$ is almost geometrically finite.
\marginwrite{algebraic lemma}
\end{proposition}

\begin{proof}{} Consider a subgroup $\Z^n$ of finite index in $V$.
Under conjugation by elements of $G$, its orbit consists of finitely
many subgroups, whose intersection is a normal subgroup isomorphic to
$\Z^n$. Let $Q'$ be the quotient of $G$ by this subgroup. Since $Q'$
maps onto $Q$ with finite kernel, it is almost geometrically finite.
Therefore we may assume that $V$ itself is free abelian of rank~$n$.
From proposition~2.2 of \cite{KLR}, there exists a Seifert
construction for this data, that is, a properly discontinuous action
of $G$ on $\R^n\times L$, where $L$ is any contractible simplicial
complex on which $Q$ acts properly discontinuously. The action takes
$\R^n$-fibers (over points of $L$) to $\R^n$-fibers, and the quotient
of each $\R^n$-fiber is a quotient of a compact flat manifold by a
finite group action. Therefore $(\R^n\times L)/G$ is compact, so $G$
is almost geometrically finite.
\end{proof}

\section [Laudenbach's Theorem] {Laudenbach's Theorem}
\label{lauden}\marginwrite{lauden}

Throughout this section let $M$ be a Haken 3-manifold, $M\neq D^3$,
and let $G\neq S^2$ be a compact connected properly imbedded
incompressible surface in $M$. Let $m_0$ be a basepoint in $G$. If $G$
is not closed, choose $m_0$ in $\partial G$. Denote by $\Imb^0(G,M)$
the space of smooth imbeddings of $G$ in $M$ that take $m_0$ to $m_0$
and $\partial G$ to $\partial G$. The inclusion map of $G$ into $M$ is
understood to be the basepoint of $\Imb^0(G,M)$. The following result
is proven on pp.~49-62 of~\cite{Lau}.

\begin{theorem}{imbeddings}{} $\pi_1(\Imb^0(G,M))\eq\set{1}$.
\marginwrite{imbeddings}
\end{theorem}

\noindent From this we deduce the result we will need.

\begin{theorem}{Laudenbach}{}
Assume that $G$ is not an annulus in a solid torus or a
boun\-dary-parallel disc. Let $f$ and $g$ be diffeomorphisms of $M$ which
preserve $G$ and fix $m_0$. Let $H$ be an isotopy from $f$ to $g$,
whose trace at $m_0$ lies in $\pi_1(G,m_0)$. Then $H$ is deformable
relative to $M\times \partial I$ to an isotopy through diffeomorphisms
that preserve~$G$. Moreover,
\begin{enumerate}
\item[{\rm(i)}]
If $H$ is relative to $\partial G$, then the new isotopy may be chosen
to be relative to $\partial G$.
\item[{\rm(ii)}]
If $H$ is relative to $\partial M$, then the new isotopy may be chosen
to be relative to $\partial M$.
\item[{\rm(iii)}]
If $H$ is relative to $\partial M$, and $f$ and $g$ agree on $G$, then
the new isotopy may be chosen to be relative to $G\cup\partial M$.
\item[{\rm(iv)}]
If $F$ is an incompressible surface (not necessarily connected) in
$M$, disjoint from $G$, such that $H$ preserves $F$, then the new
isotopy may be chosen to agree with $H$ on~$G$.
\end{enumerate}
\marginwrite{Laudenbach}
\end{theorem}

\noindent Note that by virtue of part (iv), the theorem applies to $G$
that are not connected, provided that the trace condition is satisfied
for a basepoint in each component of $G$.

\begin{proof}{}
Replacing $f$ by $g^{-1}f$, we may assume that $g$ is the identity. If
$\partial G\neq\emptyset$ and the hypothesis of one of cases (i),
(ii), or (iii) holds, then $H$ is already an isotopy relative to
$m_0$. Otherwise, we may change $f$ by isotopy in a neighborhood of
$G$ so that the trace of $H$ at $m_0$ is trivial, and then
(using~\cite{McCarty}) so that $H$ is an isotopy relative to $m_0$.
Consequently $H$ induces the identity automorphism on $\pi_1(M,m_0)$.

Let $f_1$ be the restriction of $f$ to $G$. We will first reduce to
the situation when $f_1$ is the identity. In case (iii), this already
holds. In cases (i) and (ii), $H$ is an isotopy relative to $\partial
G$ or $\partial M$. As in lemma~7.3 of~\cite{Wald}, $f_1$ is homotopic
to the identity preserving $\partial G$, so $f_1$ is isotopic to the
identity. Choose a basepoint in each component of $\partial G$.
Suppose that $\gamma$ is an arc in $G$ connecting two of them. The
restriction of $H$ to $\gamma$ gives an isotopy in $M$ from
$f(\gamma)$ to $\gamma$, and since $\pi_2(M,G)$ is trivial (because
$\pi_2(M)\eq 0$ and $\pi_1(G)\to \pi_1(M)$ is injective), $f(\gamma)$
is homotopic and hence isotopic to $\gamma$ in~$G$. Thus by changing
$f$ by isotopy preserving $G$ and relative to $\partial M$, we may
assume that $f_1$ preserves a set of arcs (disjoint except for their
endpoints) connecting the basepoints in the boundary circles. This,
together with the fact that $f_1$ is isotopic to the identity relative
to $m_0$, implies that $f_1$ is isotopic to the identity relative to
$\partial G$. After changing $f$ by isotopy preserving $G$ and
relative to $\partial M$, we may assume that $f_1$ is the identity.
Finally, suppose that none of (i), (ii), or (iii) holds. The trace
condition implies that $f$ induces the identity automorphism on
$\pi_1(M,m_0)$. Since $f(G)\eq G$, the restriction $f_1$ of $f$ to $G$
induces the identity automorphism on $\pi_1(G,m_0)$. Therefore $f_1$
is homotopic to the identity. Any orientation-preserving
diffeomorphism of a compact 2-manifold that is homotopic to the
identity is isotopic to the identity, unless $G$ is a disc or annulus
and the diffeomorphism is orientation-reversing. If $f_1$ is
orientation-reversing, then $f$ must reverse the sides of $G$. If $G$
were a disc, then since $f$ induces the identity automorphism on
$\pi_1(M,m_0)$, $M$ would have to be a 3-ball, so $G$ would be a
boundary-parallel disc. If $G$ were an annulus, then for the same
reason $M$ would have to be a solid torus, also excluded. So $f_1$ is
isotopic to the identity, and since $f_1$ induces the identity
automorphism on $\pi_1(G,m_0)$, we may assume that the isotopy
preserves $m_0$. Extending such an isotopy of $f_1$ to an isotopy of
$f$, we may assume that $f_1$ is the identity on~$G$.

Next we need to achieve that $H$ preserves $\partial G$. In (i), (ii),
or (iii), no alteration is needed. Otherwise, let $C$ be the component
of $\partial G$ that contains $m_0$, let $B$ be the component of
$\partial M$ that contains $C$, and consider the trace $\tau$ of
$H\vert_B$ at $m_0$. Note that $C$ is not contractible in $B$, since
$G$ is not a boundary-parallel disc. Suppose that $\tau$ is not
trivial in $\pi_1(B)$. Suppose for contradiction that $\tau$ does not
lie in $\pi_1(C)$. Since $f$ is the identity on $C$, $\tau$ lies in
the normalizer of $\pi_1(C)$ in $\pi_1(B)$. Therefore $B$ is a torus,
and since $\tau$ is essential in $\pi_1(B)$ and trivial in $\pi_1(M)$,
$M$ must be a solid torus. Since $\tau$ is not in $\pi_1(C)$, $C$ is
not contractible in $\pi_1(M)$, so $G$ is an annulus, the case
excluded by hypothesis. So we may assume that $\tau$ lies in
$\pi_1(C)$. Since $\tau$ is trivial in $\pi_1(M)$, $G$ must be a disc.
There is an isotopy of $M$ that preserves $G$ and moves $m_0$ around
$C$. Juxtaposing $H$ with the correct multiple of this isotopy, we may
assume that $\tau$ is trivial in $\pi_1(B)$. By~\cite{Gramain},
$\pi_1(\Imb((C,m_0),(B,m_0)))$ is trivial unless $B$ is a 2-sphere, a
case excluded since $G$ is not a boundary-parallel disc, so $H$ may be
deformed in a neighborhood of $B$ so that it preserves $C$ at each
level. Repeating for the other components of $\partial G$, we may
assume that $\partial G$ is preserved at each level of~$H$.

The restriction $j_t$ of $H$ to $G\times I$ defines a loop in
$\pi_1(\Imb^0(G,M))$. Applying theorem~\ref{imbeddings}, this loop is
contractible, so there exists a 2-parameter family $j_{t,s}$, $0\leq
s,t\leq 1$, such that $j_{t,0}\eq j_t$, and $j_{t,1}$, $j_{0,s}$, and
$j_{1,s}$ are the inclusions for each $t$ and $s$. Define $J_{t,0}\eq
H_t$, $J_{0,s}\eq f$, and $J_{1,s}\eq 1_M$. By the isotopy extension
theorem (i.~e.~the fact that $\Diff^0(M)\to \Imb^0(G,M)$ is a Serre
fibration) this extends to a 2-parameter family $J_{t,s}$ of
diffeomorphisms of $M$. Letting $K_t\eq J_{t,1}$, we have an isotopy
from $f$ to $1_M$ relative to $G$, and the existence of an isotopy
preserving $G$ together with statements (i), (ii), and (iii) are
established.

For (iv), notice that all of our alterations to $H$ may be performed
so as not to change $H$ outside a neighborhood of $G$. Therefore if
$F$ is another incompressible surface (not necessarily connected)
which is preserved by $H$, we cut along $F$ and apply the previous
argument to obtain a new isotopy agreeing with $H$ on~$F$.
\end{proof}

There is a 2-dimensional analogue of theorem~\ref{Laudenbach}.

\begin{theorem}{Laudenbach2}{}
Let $G\neq S^2$ be a surface and $k$ an arc or circle essentially
imbedded in $G$. If $k$ is a circle, let $m_0$ be a basepoint in $k$.
Let $f$ and $g$ be diffeomorphisms of $G$ which preserve $G$. Let $H$
be an isotopy from $f$ to $g$ through diffeomorphisms fixing $m_0$, if
$k$ is a circle, or preserving $\partial k$, if $k$ is an arc. Then
$H$ is deformable relative to $G\times \partial I$ to an isotopy
through diffeomorphisms that preserve~$k$. Moreover,
\begin{enumerate}
\item[{\rm(i)}]
If $H$ is relative to $\partial M$, then the new isotopy may be chosen
to be relative to $\partial M$.
\item[{\rm(ii)}]
If $H$ is relative to $\partial M$, and $f$ and $g$ agree on $G$, then
the new isotopy may be chosen to be relative to $G\cup\partial M$.
\item[{\rm(iii)}]
If $\ell$ is a $1$-manifold in $G$, disjoint from $k$, such that $H$
preserves $\ell$, then the new isotopy may be chosen to agree with $H$
on~$\ell$.
\end{enumerate}
\marginwrite{Laudenbach2}
\end{theorem}

\noindent The proof is analogous to the proof of
theorem~\ref{Laudenbach}, but much simpler.

\section [Exceptional cases] {Exceptional Seifert-fibered 3-manifolds}
\label{except}\marginwrite{except}

In the next section we will give a general treatment of the mapping
class groups of Seifert-fibered Haken 3-manifolds, but there are a few
exceptional cases to which it will not apply. We address those cases
in the present section. The manifolds (all assumed to be admissibly
fibered with complete boundary pattern) are:

\begin{enumerate}
\item[{(E1)}] The $S^1$-bundle over the annulus, with boundary pattern
$\overline{\uu{\emptyset}}$.
\item[{(E2)}] The $S^1$-bundle over the M\"obius band, with boundary
pattern $\overline{\uu{\emptyset}}$.
\item[{(E3)}] An $S^1$-bundle over the torus.
\item[{(E4)}] An $S^1$-bundle over the Klein bottle.
\item[{(E5)}] The Hantzsche-Wendt manifold, which is
the closed flat 3-manifold with Seifert invariants
$\{\,-1;\,(n_2,1);\,(2,1),\,(2,1)\,\}$ (see \cite{Orlik}~pp.~133,~138,
\cite{C-V} pp.~478-481, \cite{H-W,W1,Zimmermann}).
\item[{(E6)}] A Haken manifold which fibers over $S^2$ with three
exceptional orbits.
\end{enumerate}

\begin{proposition}{bundle exceptions}{} Let $(\Sigma,\uu{\sigma})$ be
admissibly fibered as a Seifert 3-manifold of one of the exceptional
types (E1)-(E4). Then ${\cal H}(\Sigma,\uu{\sigma})$ is almost
geometrically finite.
\marginwrite{bundle exceptions}
\end{proposition}

\begin{proof}{} For (E1), from proposition 3.4.1 of \cite{jdg} we have
${\cal H}(\Sigma,\uu{\sigma})$ isomorphic to $\Z/2\times \GL(2,\Z)$
(the $\Z/2$ factor is generated by reflection in the $I$-fibers of the
$I$-bundle structure). So ${\cal H}(\Sigma,\uu{\sigma})$ is virtually
free. For (E2), ${\cal H}(\Sigma,\uu{\sigma})$ is finite, by
proposition 3.4.2 of \cite{jdg}. For (E3), if the Euler class is zero
then $\Sigma$ is the 3-torus, with ${\cal
H}(\Sigma)\cong\Out(\pi_1(\Sigma))\cong\GL(3,\Z)$. If the Euler class
is $n$, then $\pi_1(\Sigma)\cong \langle x,y,t\vbar
[x,t]=[y,t]=1,[x,y]=t^n\rangle$. The center is $\Z$ generated by $t$,
and the quotient of $\pi_1(\Sigma)$ by its center is $\Z\times\Z$
generated by the images of $x$ and $y$. Sending $\Out(\pi_1(\Sigma))$
to $\Aut(\Z\times\Z)=\GL(2,\Z)$ is surjective and splits (the three
automorphisms determined by sending (1)~$x$ to $xy$, $y$ to $y$, and
$t$ to $t$, (2)~$x$ to $x$, $y$ to $xy$, and $t$ to $t$, and (3)~$x$
to $y$, $y$ to $x$, and $t$ to $t^{-1}$ define the splitting).
Elements of the kernel are the automorphisms $\phi(i,j)$ that send $x$
to $xt^i$, $y$ to $yt^j$, and $t$ to $t$. Conjugation by $x$ equals
$\phi(0,n)$ and conjugation by $y$ equals $\phi(-n,0)$, so the kernel
of $\Out(\pi_1(\Sigma))\to\GL(2,\Z)$ is $\Z/n\times\Z/n$, showing that
$\hscr(\Sigma)\cong\Out(\pi_1(\Sigma))$ is almost geometrically
finite. The manifolds of type (E4) are analyzed in proposition~3.4.4
of \cite{jdg}. If the Euler class is zero then there is a homomorphism
from $\Out(\pi_1(M))$ to $\hbox{PGL}(2,\Z)$ with finite kernel. Since
$\PGL(2,\Z)$ is virtually free, it is almost geometrically finite. If
the Euler class is nonzero, $\Out(\pi_1(\Sigma))$ is finite.
\end{proof}

\begin{proposition}{other exceptions}{} Let $M$ be one of the
exceptional types (E5)-(E6). Then ${\cal H}(M)$ is finite.
\marginwrite{bundle exceptions}
\end{proposition}

\begin{proof}{} For (E5), it is from \cite{C-V} (although the correct
structure of the group was later given in \cite{Zimmermann}). The
lemma in section~3.4 of \cite{jdg} gives case~(E6).
\end{proof}

\section [Fibered 3-manifolds] {Fibered 3-manifolds}
\label{seifert}\marginwrite{seifert}

We first treat the case of $I$-bundles, over surfaces which are not
necessarily connected.

\begin{lemma}{ibundle mapping class group}{} Suppose that
$(\Sigma,\s)$ is admissibly $I$-fibered over the compact 2-manifold
$(B,\b)$. Let $p\colon\Sigma\to B$ be the projection. Then
$\hscr(\Sigma,\uu{\sigma})$ is isomorphic to a semidirect product
$F\circ\hscr(B,\b)$ where $F$ is a direct sum of copies of
$\Z/2$, one for each component of $\Sigma$, and the action of an
element of $\hscr(B,\b)$ on $F$ is to permute the copies $\Z/2$
exactly as it permutes the corresponding components of~$\Sigma$.
\marginwrite{ibundle mapping class group}
\end{lemma}

\begin{proof}{} Let $i\colon (B,\b)\to(\Sigma,\s)$ be the $0$-section
of the $I$-bundle (where $I$ is regarded as $[-1,1]$ and the structure
group is reduced to $\Z/2$ generated by reflection in $I$, and the
twisting is given by the orientation homomorphism). For each
component of $\Sigma$ there is the involution given by reflection in
the $I$-fibers, and this is not admissibly isotopic to the identity
since it is orientation-reversing. These give the generators of~$F$.

Define $j\colon\hscr(B,\b)\to\hscr(\Sigma,\s)$ by extending
diffeomorphisms on $i(B)$ to $\Sigma$ linearly in each fiber, choosing
the unique way to do this that is orientation-preserving on $\Sigma$.
To see that this is injective, let $\langle h\rangle$ be an element of
$\hscr(B,\b)$ such that $j(\langle h\rangle)$ is trivial in
$\hscr(\Sigma,\s)$. Then $h$ is orientation-preserving, since
otherwise the restriction of $j(h)$ of the lid or lids of $\Sigma$
cannot be isotopic to the identity. Projecting an admissible isotopy
from $j(h)$ to $1_\Sigma$ down to $i(B)$ gives an admissible homotopy
from $h$ to $1_B$. This implies that $h$ is admissibly isotopic to
$1_B$ (see for example lemma~2.19 of~\cite{C-M}). Corollary~5.9
of~\cite{Joh} shows that the image of $j$ is the entire group of
orientation-preserving mapping elements of $\hscr(\Sigma,\s)$. Then,
it is clear that the subgroups $F$ and $j(\hscr(\Sigma,\s))$ generate
$\hscr(\Sigma,\s)$, and $F$ is normal, and the lemma follows.
\end{proof}

Throughout the remainder of this section, $(\Sigma,\uu{\sigma})$ will
denote a Haken 3-manifold with complete and useful boundary pattern,
which admits an admissible Seifert fibering over $(B,\uu{b})$. A
diffeomorphism of $\Sigma$ is called {\it fiber-preserving} if it
carries each fiber $\Sigma$ to a fiber of $\Sigma$, and is called
{\it vertical} if it takes each fiber to itself. By
$\hscr_f(\Sigma,\uu{\sigma})$ we indicate the mapping classes of
fiber-preserving diffeomorphisms (that is, fiber-preserving
diffeomorphisms modulo isotopy through fiber-preserving
diffeomorphisms). There is a natural homomorphism
$\hscr_f(\Sigma,\s)\to\hscr(\Sigma,\s)$.

\begin{theorem}{fiber-preserving}{} If $(\Sigma,\s)$ is not one of
(E1)-(E6), then $\hscr_f(\Sigma,\s)\to\hscr(\Sigma,\s)$ is an
isomorphism.
\marginwrite{fiber-preserving}
\end{theorem}

\begin{proof}{} Provided that $(\Sigma,\s)$ is not an exception
(E1)-(E5), \cite{W1} or theorem~8.1.7 of \cite{Orlik} shows that the
homomorphism is surjective. If it is not an exception (E6), then it has
a vertical incompressible surface, and the argument of p.~85-86 of
\cite{Wald} shows that it is injective.
\end{proof}

Define $\hscr^0(\Sigma,\s)$ to be the elements of $\hscr(\Sigma,\s)$
that contain a vertical diffeomorphism. As in lemma~25.2 of \cite{Joh}
(see also lemma~3.5.3 of~\cite{jdg}), we have the following
calculation.

\begin{proposition}{vertical mcg}{} Let $\Sigma$ be Seifert-fibered
over $(B,\uu{b})$, with no component of $\Sigma$ an exceptional case
(E1)-(E6). Then $\hscr^0_+(\Sigma,\s)\cong H_1(B,\abs{\uu{b}})$.
\marginwrite{vertical mcg}
\end{proposition}

\noindent The rough idea behind this result is that the generators of
$\hscr^0_+(\Sigma,\s)$ are Dehn twists (see \cite{Joh} or section~3.3
of \cite{jdg} for a definition of Dehn twist) about vertical tori and
annuli, which obey the same homological relations as their image
circles and arcs in $B$ (in the nonorientable case there is another
type of generator supported in a neighborhood of a vertical Klein
bottle; its square is a Dehn twist about the boundary of a regular
neighborhood of the Klein bottle.

Let $E$ be the exceptional points of $B$, that is, the images of the
exceptional orbits of $\Sigma$. This is a finite subset of the
interior of $B$. Denote by $\rho\colon
\hscr(\Sigma,\s)\to\hscr(B-E,\b)$ the homomorphism induced by
projection of fiber-preserving homomorphisms to the base surface.
Define $\hscr_0(B-E,\uu{b})$ to be the subgroup of $\hscr(B-E,\uu{b})$
consisting of the classes $\langle f\rangle$ such that $f$ is
admissibly isotopic to a map which is the identity on
$\abs{\uu{b}}\eq\partial B$, and $f$ permutes the punctures of $B$
trivially. Since $\uu{b}$ consists of arcs and circles (and since
$B$ is of finite type), $\hscr_0(B-E,\uu{b})$ has finite index
in~$\hscr(B-E,\uu{b})$.

\begin{proposition}{homeos of base}{}
The image of $\rho\colon \hscr_f(\Sigma,\s)\to\hscr(B-E,\b)$ contains
$\hscr_0(B-E,\uu{b})$, hence has finite index in $\hscr(B-E,\uu{b})$.
If $\partial B\neq\emptyset$, then there exists a homomorphism
$s\colon \hscr_0(B-E,\uu{b})\to \hscr_f(\Sigma,\s)$ such that $\rho s$
is the identity, and such that each $s(\langle f\rangle)$ has a
representative which is the identity on $\partial \Sigma$.
\marginwrite{homeos of base}
\end{proposition}

\begin{proof}{} The proposition follows from the special case
when $B$ is connected. The connected case appears a proposition~25.3 of
\cite{Joh} and theorem~3.5.2 of \cite{jdg}. The splitting is
constructed using a cross-section $i\colon B-E\to
\Sigma-\widetilde{E}$, where $\widetilde{E}$ is the union of the
exceptional fibers. For $\langle f\rangle\in\hscr_0(B-E,\uu{b})$, one
may assume that $f$ is the identity on $\partial B$, so
(since $s(f)$ is always selected to be orientation-preserving) that
$s(f)$ is the identity on $i(\partial B)$ and hence on~$\partial\Sigma$.
\end{proof}

\noindent Propositions~\ref{vertical mcg} and~\ref{homeos of base},
combined with lemma~\ref{2-manifold mcg}, yield immediately the
following:

\begin{theorem}{Sfs mcg}{} Let
$(\Sigma,\uu{\sigma})$ be an admissibly Seifert-fibered Haken
3-manifold with complete and useful boundary pattern $\sigma$. Assume
that $(\Sigma,\uu{\sigma})$ is not one of the exceptional manifolds
(E1)-(E6). Then there is an exact sequence
$$1\longrightarrow V\longrightarrow
{\cal H}(\Sigma,\uu{\sigma})\longrightarrow Q\longrightarrow1$$
\noindent where $V$ has a finitely generated abelian subgroup of
finite index, and $Q$ is almost geometrically finite.
\marginwrite{Sfs mcg}
\end{theorem}

\noindent Applying proposition~\ref{algebraic lemma}, we have

\begin{corollary}{Sfs mcg case}{} Let
$(\Sigma,\uu{\sigma})$ be an admissibly Seifert-fibered Haken
3-manifold with complete and useful boundary pattern $\sigma$. Assume
that $(\Sigma,\uu{\sigma})$ is not one of the exceptional manifolds
(E1)-(E6). Then ${\cal H}(\Sigma,\uu{\sigma})$ is almost geometrically
finite.
\marginwrite{Sfs mcg case}
\end{corollary}

\section [Haken manifolds] {Haken manifolds}
\label{haken}\marginwrite{haken}

Throughout this section we assume that $(M,\m)$ is a Haken manifold
with a complete and useful boundary pattern. We allow the possibility
that $\partial M$ is empty. Denote by $\Sigma$ the characteristic
submanifold of $(M,\m)$. Let $\s$ be the proper boundary pattern on
$\Sigma$. The following result was proved in~\cite{jdg}, but we
present a less abbreviated proof here. A reference for the $\hbox{\it
Sol}$ geometry is theorem~5.5 of~\cite{Scottgeom}).

\begin{proposition}{preserving sigma mcg}{} Suppose that $M$ is not a
torus bundle over the circle that admits a $\hbox{\it Sol}$ structure.
Then $\hscr(M,\Sigma,\m)\to \hscr(M,\m)$ is an isomorphism.
\marginwrite{preserving sigma mcg}
\end{proposition}

\begin{proof}{} Since the characteristic submanifold is unique up to
admissible isotopy, the homomorphism is surjective. For injectivity,
let $\langle f \rangle\in \hscr(M,\Sigma,\m)$ and suppose that
$H\colon M\times I\to M$ is an admissible isotopy from $f$ to $1_M$.
We must find an isotopy that preserves the frontier of~$\Sigma$.

Let $F$ be a component of the frontier of $\Sigma$. We claim that
$f(F)\eq F$. Suppose not. The restriction of $H$ to $F\times I$ is an
admissible map from an $I$-bundle or Seifert fiber space into
$(M,\m)$, so by proposition~13.1 of \cite{Joh},
it is admissibly homotopic into $\Sigma$.
That is, two components of the frontier of
$\Sigma$ are admissibly homotopic in $\Sigma$.
By proposition~19.1 of \cite{Joh}, these components must
be admissibly parallel in $\Sigma$,
that is, the component of $\Sigma$ containing $F$ is of the form
$F\times I$ and $f$ interchanges of the components of its boundary.
Therefore $f$ is isotopic to a diffeomorphism that preserves $F$ and
interchanges its sides. Since $f$ is isotopic to the identity,
lemma~7.4 of~\cite{Wald} shows that this is impossible.

By isotopy preserving $\Sigma$, we may assume that $f$ fixes a
basepoint $m_0$ in the interior of~$F$. We claim that the trace of $H$
at $m_0$ is in the subgroup $\pi_1(F,m_0)$. When $\partial M$ is
nonempty, corollary~18.2 of~\cite{Joh} applies to prove the claim.
When $\partial M$ is empty, the argument in lemma~18.1 of~\cite{Joh}
shows that if the claim is false then the components of $\Sigma$ and
$\overline{M-\Sigma}$ adjacent to $F$ are each diffeomorphic to the product
of the torus and an interval. By maximality of $\Sigma$, this is only
possible when $M$ is a torus bundle over the circle which admits a
{\it Sol} structure, which is excluded by hypothesis. This establishes
the claim.

Since $F$ is a square, annulus, or torus, there is an isotopy on $F$
from the identity to the identity whose trace is equal to the trace of
$H$. So it is possible to change $f$ by an admissible isotopy with
support in a neighborhood of $F$, so that the trace of the isotopy
from $f$ to the identity of $M$ is trivial. Then $f$ must induces the
identity automorphism on $\pi_1(M,m_0)$. Since $f$ preserves $F$, it
also induces the identity isomorphism on $\pi_1(F,m_0)$.

Let $f_1$ be the restriction of $f$ to $F$. We
will show that $f_1$ preserves each element of the boundary
pattern of $F$. If $F$ is a torus, there is nothing to prove. Suppose
$F$ is an annulus, so that the boundary pattern consists of the two
boundary circles. If $f_1$ interchanges them, then since $f_1$ induces the
identity automorphism it must be orientation-reversing. But
$f$ preserves the sides of $F$ (because $f$ preserves
$\Sigma$, or alternatively using lemma~7.4 of~\cite{Wald} again), so $f$
would be orientation-reversing and could not be isotopic to the
identity. Suppose that $F$ is a square. Its boundary pattern consists
of the four edges. Suppose for contradiction that $f_1$ moves some
edge to a different edge. Since $f$ is admissibly homotopic to the
identity, it must preserve each element of $\m$. Since adjacent edges
cannot lie in the same element of $\m$, $f_1$ must interchange a pair
of opposite edges. If it interchanges one pair of opposite edges, but
preserves each of the other two edges, then $f_1$ is
orientation-reversing, a contradiction as in the case when $F$ is an
annulus. Therefore $f$ must interchange both pairs of opposite edges.
Since $F$ is a square, the component of $\Sigma$ containing $F$ must
be an $I$-bundle, and since opposite edges of $F$ lie in the same
component of the boundary pattern, one pair lies in the lid and the
other pair are joined by a square $S$ which is contained in an element
$F'$ of $\m$. Now the $I$-bundle cannot be fibered over the disc,
because since $F$ and $S$ are sides which meet in two fibers, the
$I$-bundle would be fibered over a 2-faced disc and would not be
essentially imbedded. This is impossible since it is a component of
$\Sigma$. So the center circle of the annulus $F\cup S$ is essential
in $M$. Since $f$ interchanges opposite edges of $F$, it must send
this element of $\pi_1(M)$ to its negative, which is impossible since
$f$ induces the identity automorphism on~$\pi_1(M,m_0)$.

Fix a component $F$ of the frontier of $\Sigma$. Suppose $F$ is a
torus. Since the trace of the isotopy at $m_0$ is trivial,
theorem~\ref{Laudenbach} applies to show that the isotopy from $f$ to
the identity may be deformed to preserve $F$ at each stage.

Suppose $F$ is an annulus. Consider a basepoint $b_0$ in a boundary
circle $C$ of $F$. Since $f_1$ is the identity on $F$, and the trace of
$H$ at $m_0$ is trivial, the trace of $H$ at $b_0$ is also trivial in
$\pi_1(M,b_0)$. Since the elements of $\m$ are incompressible, the
trace at $b_0$ is trivial in $\pi_1(G)$ where $G$ is the element of
$\m$ that contains $b_0$. Applying theorem~\ref{Laudenbach2}, we may
deform the isotopy admissibly in a neighborhood of $\partial M$ so
that $C$ preserved at each level of the isotopy. Repeating, we assume
that all of $\partial F$ is preserved, and then apply
theorem~\ref{Laudenbach}.

Now suppose $F$ is a square. Let $b_0$ be a corner where two edges
meet, so $b_0\in G_1\cap G_2$ for two elements of the boundary
pattern. Again, the trace of $H$ at $b_0$ is trivial in $\pi_1(G_1)$
and $\pi_1(G_2)$. If it is not trivial in the fundamental group of the
component $k$ of $G_1\cap G_2$ that contains it, then that component is a
circle and $G_1$ and $G_2$ must be discs. Then $M$ would be a 3-ball
with boundary pattern $\set{G_1,G_2}$ and $M$ could not contain the
essentially imbedded square $F$. We conclude that the trace of $H$ at
$b_0$ is trivial in $\pi_1(k)$, so we may assume that $H$ preserves
$b_0$. Repeating, we assume that $H$ preserves each of the four
corners of $F$. Now using theorem~\ref{Laudenbach2} we may assume that
$H$ preserves each of the four edges of $F$, and apply
theorem~\ref{Laudenbach} to assume that $H$ preserves~$F$.

Repeating this for each component of the frontier of $\Sigma$, not
disturbing those already adjusted, shows that $f$ was the trivial
element of~$\hscr(M,\Sigma,\m)$.
\end{proof}

From now until when we reach theorem~\ref{main theorem}, we assume
that $M\neq \Sigma$. Let $H\eq \ol{M-\Sigma}$, and let $\uu{h}$ be the
proper boundary pattern on~$H$. Define $R_H$ to be the image of the
restriction $\hscr(M,\Sigma,\m)\to \hscr(H,\uu{h})$. From lemma~4.2.1
of~\cite{jdg} we have the following fact.

\begin{lemma}{finite image in mcg(H)}{} $R_H$ is finite.
\marginwrite{finite image in mcg(H)}
\end{lemma}

\noindent Actually, according to theorem~\ref{FMCGT},
$\hscr(H_i,\uu{h_i})$ is itself finite for all components of
$(H,\uu{h})$ except the case when $(H_i,\uu{h_i})\eq (T^2\times
I,\overline{\uu{\emptyset}})$. For these a special argument is given
in \cite{jdg}. The idea is that on the adjacent Seifert-fibered
component(s) of $\Sigma$, any diffeomorphism must be fiber-preserving
up to isotopy. The fibers are linearly independent in $\pi_1(T^2\times
I)\cong\Z\times\Z$, and only $\pm I$ can preserve two linearly
independent elements of $\Z\times\Z$, allowing at most two
possibilities for the restriction to~$\hscr(H_i,\uu{h_i})$.

Since $M\neq\Sigma$, the components of $\Sigma$, with their proper
boundary patterns, are of four kinds.
\begin{enumerate}
\item[{\rm (i)}]
The components which can be admissibly fibered as
$I$-bundles with their lids in $\partial M$. Their union will be
denoted by~$I$.
\item[{\rm (ii)}]
Those that are diffeomorphic to $S^1\times S^1\times I$ with boundary
pattern $\ol{\uu{\phi}}$, but are not as in (i) (that is, exceptional
case (E1)). These must be Seifert-fibered (since as $I$-bundles they
could not have both lids in $\partial M$), and must have frontier
equal to one or both of their boundary components. Their union will be
denoted by~$T$.
\item[{\rm (iii)}]
Those that are diffeomorphic to the twisted $I$-bundle over the Klein
bottle with boundary pattern $\ol{\uu{\phi}}$ (that is, exceptional
case (E2)). Their union will be
denoted by~$K$.
\item[{\rm (iv)}]
The components that are not as in (i), (ii), or (iii). These are
Seifert-fibered, since they are not as in (i), and their fiberings are
unique up to admissible isotopy, since they are not as in (ii)
or~(iii). Their union will be denoted by~$S$.
\end{enumerate}

The proper boundary pattern of $I$ is denoted $\uu{i}$, and similarly
for $K$, $T$, and $S$. Clearly each of $I$, $K$, $T$, and $S$ is
preserved by each element of $\hscr(M,\Sigma,\m)$. Let $R_T$ be the
image of the restriction $\hscr(M,\Sigma,\m)\to \hscr(T,\uu{t})$.

\begin{lemma}{image in hscr(T) is finite}{} $R_T$ is finite.
\marginwrite{image in hscr(T) is finite}
\end{lemma}

\begin{proof}{} Consider the composition
$$\hscr(M,\Sigma,\m)\to R_T\to \hscr(\Fr(T))\ ,$$

\noindent where $\Fr(T)$ is the frontier.
Since this equals the composition
$$\hscr(M,\Sigma,\m)\to \hscr(H,\uu{s})\to
\hscr(\Fr(T))\ ,$$

\noindent and $\hscr(M,\Sigma,\m)\to \hscr(H,\uu{s})$ has
finite image by lemma~\ref{finite image in mcg(H)}, it follows
that $R_T\to \hscr(\Fr(T))$ has finite image. But it is also injective,
since any diffeomorphism of $T$ isotopic to the identity on a
boundary component is isotopic to the identity. Therefore $R_T$ is
finite.
\end{proof}

\begin{lemma}{hscr(K) is finite}{} $\hscr(K,\uu{k})$ is finite.
\marginwrite{hscr(K) is finite}
\end{lemma}

\begin{proof}{} Since $(K,\uu{k})$ fibers admissibly as an $I$-bundle,
$\uu{k}\eq\overline{\uu{\emptyset}}$. By lemma~\ref{ibundle mapping
class group}, if $X$ is the twisted $I$-bundle over the Klein bottle
$Y$, we have $\hscr(X,\overline{\uu{\emptyset}})\cong\hscr(Y)$, and
by~\cite{Lickorish} the latter is $\Z/2\times\Z/2$. Since $K$ has
finitely many components, $\hscr(K,\uu{k})$ is also finite.
\end{proof}

Consider the restriction homomorphism
$$\hscr(M,\Sigma,\m)\to R_H\times R_T\times
\hscr(K,\uu{k})\times\hscr(I,\uu{i})\times\hscr(S,\uu{s})\ .$$

\noindent Let $(B,\uu{b})$ be the quotient surface of $(S,\uu{s})$,
and let $E$ be the image of the exceptional fibers. By
theorem~\ref{fiber-preserving} and proposition~\ref{homeos of base},
there is a homomorphism $\hscr(S,\uu{s})\to\hscr(B-E,\uu{b})$, and
composing with this in the $\hscr(S,\uu{s})$ factor, we obtain a
homomorphism
$$\Psi\colon\hscr(M,\Sigma,\m)\to R_H\times R_T\times
\hscr(K,\uu{k})\times\hscr(I,\uu{i})\times\hscr(B-E,\uu{b})\ .$$

\begin{proposition}{Psi has image of finite index}{} The image of $\Psi$
has finite index.
\marginwrite{Psi has image of finite index}
\end{proposition}

\begin{proof}{} By lemmas~\ref{finite image in mcg(H)}, \ref{image in
hscr(T) is finite}, and~\ref{hscr(K) is finite}, we need only examine
$\hscr(I,\uu{i})$ and $\hscr(B-E,\uu{b})$. Let $\hscr_0(I,\uu{i})$  be
the subgroup of $\hscr(I,\uu{i})$ consisting of the elements
containing representatives which are the identity on the frontier of
$I$, and recall the subgroup $\hscr_0(B-E,\uu{b})$ of
$\hscr(B-E,\uu{b})$ defined shortly before proposition~\ref{homeos of
base}. Since the mapping class groups of the square, annulus, arc, and
circle are finite, these are finite index subgroups. It suffices to
show that $\hscr_0(I,\uu{i})\times \hscr_0(B-E,\uu{b})$ is contained
in the image of~$\Psi$.

Proposition~\ref{homeos of base} shows that the image of the
homomorphism $\rho\colon\hscr_f(S,\uu{s}) \to\hscr(B-E,\b)$ contains
$\hscr_0(B-E,\uu{b})$, and that there is a homomorphism
$s\colon\hscr_0(B-E,\uu{b})\to\hscr_f(S,\uu{s})$ with $\rho s$ equal to
the identity. Moreover, each $s(\langle h\rangle)$ can be represented
by a diffeomorphism which is the identity on $\partial S$. So given an
element in $\hscr_0(I,\uu{i})\times \hscr_0(B-E,\uu{b})$, it can be
represented on $I\cup S$ by a diffeomorphism which is the identity on
the frontier of $I\cup S$. Extending this to $M$ using the identity
map on $M-(I\cup S)$ produces an element of $\hscr(M,\Sigma,\m)$ that
$\Psi$ carries to the given element.
\end{proof}

Define $\kscr(M,\Sigma,\m)$ to be the kernel of $\Psi$.

\begin{proposition}{abelian kernel}{} $\kscr(M,\Sigma,\m)$ is
finitely generated and abelian.
\marginwrite{abelian kernel}
\end{proposition}

\begin{proof}{}
Let $g$ be a vertical diffeomorphism of $S$ which is isotopic to the
identity on $\Fr(S)$. Each such $g$ extends to a diffeomorphism of $M$
which is the identity outside a regular neighborhood of $S$. The image
of $\kscr(M,\Sigma,\m)$ in $\hscr(S,\uu{s})$ is contained in the
subgroup of $\hscr^0_+(S,\uu{s})$ consisting of elements representable
by maps which are isotopic to the identity on the frontier of $S$. By
proposition~\ref{vertical mcg}, $\hscr^0_+(S,\uu{s})$ is finitely generated
and abelian, hence so is the image of $\kscr(M,\Sigma,\m)$. Choose a
set of generators for the image, and for each an extension to $M$. Let
$X_1$ be the resulting subset of $\hscr(M,\Sigma,\m)$. These commute
with each other and with any diffeomorphism supported in a regular
neighborhood of $\Fr(H)$. Choose a set of generators $X_2$ (Dehn
twists about components of $\Fr(H)$) for the elements of
$\hscr(M,\Sigma,\m)$ representable by diffeomorphisms supported in a
regular neighborhood of $\Fr(H)$. The set $X\eq X_1\cup X_2$ is a
finite set of commuting elements, and we claim that it generates
$\kscr(M,\Sigma,\m)$. Given any element of $\kscr(M,\Sigma,\m)$,
we may compose this element with elements in $X_1$ to assume that its
restriction to $S$ is isotopic to the identity. But all elements of
$\kscr(M,\Sigma,\m)$ are isotopic to the identity when restricted to
each of $I$, $K$, $T$, and $H$, so when the restriction to $S$ is also
isotopic to the identity, the element is a product of elements in
$X_2$.
\end{proof}

\begin{theorem}{main theorem}{} Let $M$ be a Haken manifold and $\m$ a
boundary pattern on $M$ whose completion is useful. Then $\hscr(M,\m)$
is virtually torsionfree and almost geometrically finite.
\marginwrite{main theorem}
\end{theorem}

\begin{proof}{}
If $M$ is a torus bundle over $S^1$ that admits a $\hbox{\it Sol}$
structure, then by proposition~4.1.2 of \cite{jdg}, $\hscr(M)$ is
finite. Otherwise, by proposition~\ref{preserving sigma mcg}, we may
work instead with $\hscr(M,\Sigma,\m)$.

First we combine previous results to prove that $\hscr(M,\m)$ is
almost geometrically finite. If $M\eq \Sigma$, then lemma~\ref{ibundle
mapping class group} and corollary~\ref{Sfs mcg case} apply. Suppose that
$M\neq\Sigma$. By lemma~\ref{ibundle mapping class group},
$\hscr(I,\uu{i})$ is almost geometrically finite. By
proposition~\ref{Psi has image of finite index}, the image of $\Psi$
is almost geometrically finite. By proposition~\ref{abelian kernel},
the kernel of~$\Psi$ is finitely generated abelian. By
proposition~\ref{algebraic lemma}, $\hscr(M,\m)$ is almost
geometrically finite.

It was proven in \cite{jdg} that $\hscr(M,\m)$ is virtually
torsionfree, but we give here a proof that is much simpler.

If $M\eq\Sigma$, then lemma~3.5.9 of~\cite{jdg} completes the proof,
so we assume that $M\neq \Sigma$. Let $(Q,\uu{q})$ be obtained as
follows. Let $(B_1,\uu{b_1})$ be the base surface of $(I,\uu{i})$, and
let $(B_2,\uu{b_2})$ be the base surface of $(S,\uu{s})$, and $E$
the image of the exceptional fibers. Remove each
component $(X,\uu{x})$ from $(B_1,\uu{b_1})\cup (B_2,\uu{b_2})$ for
which $\hscr(X-E,\uu{x})$ is finite, and call the remainder
$(Q,\uu{q})$. Let
$$\pi\colon R_H\times R_T\times
\hscr(K,\uu{k})\times\hscr(I,\uu{i})\times\hscr(B_2-E,\uu{b_2})\longrightarrow
\hscr(Q-E,\uu{q})$$
\noindent be the projection, and note that $\pi$ has finite kernel.
Put $\Psi'\eq \pi \Psi$ and $V\eq(\Psi')^{-1}(F)$ where $F$ is the
subgroup of
$\hscr(I,\uu{i})$ from lemma~\ref{ibundle mapping class group}. Then
we have an exact sequence
$$1\to V\mapright{} \hscr(M,\m)\mapright{\Psi'}
\hscr(Q-E,\uu{q})$$

\noindent where $V$ is virtually abelian. Let $\gscr(Q-E,\uu{q})$ denote
the subgroup of $\hscr(Q-E,\uu{q})$ consisting of the elements that
preserve each element of $\uu{q}$ and each point of intersection of
two elements of $\uu{q}$, and let $\gscr_0(Q-E,\uu{q})$ be the
corresponding subgroup of $\hscr_0(Q-E,\uu{q})$. Let $\hscr'(M,\m)$
denote the subgroup $(\Psi')^{-1}(\gscr_0(Q-E,\uu{q}))$; as in
proposition~\ref{homeos of base}, there is a splitting $s\colon
\gscr_0(Q-E,\uu{q}) \to \hscr'(M,\m)$ so we have a semidirect product
$$\hscr'(M,\m)\eq V\circ \gscr_0(Q-E,\uu{q})\ .$$

\noindent
Let $W\eq H_1(B_1)\times V$, and consider the {\it holomorph} of $W$,
$W\circ\Aut(W)$, which is defined to be the semidirect product in
which $\Aut(W)$ acts naturally on $W$. We claim that $\Aut(W)$ is
virtually torsionfree. Write $W$ as $W_1\circ T$ where $W_1$ is free
abelian and $T$ is finitely generated. Since $T$ is characteristic,
there is a homomorphism $\Aut(W)\to\Aut(W_1)$. It is easily checked
that the kernel is a finite group of the form
$\Hom(W_1,T)\circ\Aut(T)$, and moreover it splits. Since
$\Aut(W_1)\cong \GL(n,\Z)$ is virtually torsionfree, it follows that
$\Aut(W)$ is virtually torsionfree. Since $W$ is also virtually
torsionfree, and $W$ is finitely generated, lemma~6.8 of~\cite{M-M}
implies that $W\circ\Aut(W)$ is virtually torsionfree.

The action in the semidirect product $V\circ \gscr_0(Q-E,\uu{q})$
together with the natural action of $\hscr_0(B_1)$ on $H_1(B_1)$
defines a homomorphism from $\gscr_0(Q-E,\uu{q})$ to $W\circ\Aut(W)$.
Since the latter is virtually torsionfree, it suffices
to show that the kernel is torsionfree. The kernel is
contained in the subgroup of $\gscr_0(Q-E,\uu{q})$ consisting of elements
that act trivially on $H_1(Q-E)$. Since we have removed all disc
components from $Q-E$, elements of the kernel preserve each component
of $Q-E$, so we may assume that $Q-E$ is connected. Let $P$ consist of
one point from each arc of $\uu{q}$ together with all points that are
intersections of two distinct arcs of $\uu{q}$. It is easy to see that
$\gscr_0(Q-E,\uu{q})$ is isomorphic to a subgroup of finite index in
$\hscr(Q-E\rel P)$. If $P$ is empty, then it is well-known that the
subgroup of $\hscr(Q-E)$ inducing the identity on homology is
torsionfree (see for example~\cite{Farkas-Kra}). Suppose $P$ is not
empty. Let $\Diff_1(Q-E)$ denote the subgroup of $\Diff(Q-E)$ consisting
of the elements which induce the identity on $H_1(Q-E)$ and preserve
each component of $\partial (Q-E)$. There is a fibration
$$\Diff_1(Q-E\rel P)\to\Diff_1(Q-E)\to\imb(P,\partial(Q-E))$$
\noindent where $\imb(P,\partial(Q-E))$ is the connected component of the
inclusion in the space of imbeddings $\Imb(P,\partial(Q-E))$. Clearly
$\imb(P,\partial (Q-E))\cong (S^1)^k$, where $k$ is the number of boundary
circles of $Q-E$ that contain points of $P$. From the homotopy exact
sequence of this fibration, we have an exact sequence
$$0\to\pi_1(\Diff_1(Q-E))\to\Z^k\to\hscr_1(Q-E\rel
P)\to\hscr_1(Q-E)\to 1\ $$
where $\hscr_1(Q-E\rel P)\eq \pi_0(\Diff_1(Q-E\rel P))$ and
$\hscr_1(Q-E)\eq \pi_0(\Diff_1(Q-E))$. Now $\pi_1(\Diff_1(Q-E))$ is
nontrivial only when $Q-E$ is an annulus and $k\eq 2$ (a disc, a
M\"obius band, or an annulus with $k\eq1$ have already been excluded
by the definition of $Q-E$), and in this case the generator of
$\pi_1(\Diff_1(Q-E))$ is carried to an element of the form
$(\pm1,\pm1)\in\Z\times\Z$. Therefore we have an exact sequence
$$0\to\Z^\ell\to \hscr_1(Q-E\rel P)\to\hscr_1(Q-E)\to 1\ .$$
\noindent We have already seen that $\hscr_1(Q-E)$ is torsionfree, and
therefore $\hscr_1(Q-E\rel P)$ is torsionfree.
\end{proof}

To deduce some corollaries, it is convenient to introduce a general
technique for extending results about $\hscr(M,\m)$ to relative
mapping class groups. Assume that $\m$ is complete. Suppose
$\mone\subset\m$. Take a fine triangulation of $\abs{\mone}$, which
includes as vertices the points of $\partial\abs{\uu{m_1}}$ that are
intersections of three elements of $\uu{m}$, and let $\mtwo$ be the
complete boundary pattern on $M$ consisting of the elements of
$\m-\mone$ and the 2-cells dual to the triangulation of $\abs{\mone}$.
Choose the triangulation so that at every vertex in $\abs{\mone}$, at
least four triangles meet (for example, take the barycentric
subdivision of the triangulation of $\abs{\uu{m_1}}$ without
introducing the barycenters of 1-simplices that lie in the boundary.
Then, each dual 2-cell has at least four sides. This ensures the
following.
\begin{enumerate}
\item[{\rm (i)}]
If $\m$ is useful, then $\mtwo$ is useful.
\item[{\rm (ii)}]
If $\m-\mone$ consists of the components of $\ol{\partial
M-\abs{\mone}}$, and these components are incompressible, then $\mtwo$
is useful.
\end{enumerate}

\noindent The important property of $\mtwo$ is the following.

\begin{proposition}{relative trick}{}
$\hscr(M,\m\rel\abs{\mone})$ is isomorphic to a subgroup of finite
index in $\hscr(M,\mtwo)$.
\marginwrite{relative trick}
\end{proposition}

\begin{proof}{}
Consider the natural homomorphism from $\hscr(M,\m\rel \abs{\mone})$
to $\hscr(M,\mtwo)$. We claim this is injective and has image of
finite index. Suppose $f$ and $g$ are the identity on $\abs{\mone}$
and are equivalent in $\hscr(M,\mtwo)$. Any admissible isotopy must
preserve the intersection of any collection of dual 2-cells, so must
preserve the points where three dual 2-cells in $W$ meet. Therefore we
may assume the isotopy is relative to the dual 0-cells. Also, it must
preserve the intersection of any two dual 2-cells, i.~e.~the dual
1-cells. By the Alexander trick applied to each $\hbox{dual
1-cell}\times I$, we may make the admissible isotopy to be relative to
the intersections of the dual cells. Then, by the Alexander trick
applied to each $\hbox{dual 2-cell}\times I$, we may make the isotopy
relative to all of~$\abs{\mone}$. To show the image has finite index,
let $\hscr(M,\mtwo)$ act on the union of the set of all the dual
cells. This defines a homomorphism to a finite permutation group, and
by an argument similar to the proof of injectivity one shows that any
element in the kernel is admissibly (for $\mtwo$) isotopic to the
identity on $\abs{\mone}$, so is in the image.
\end{proof}

Using $\mtwo$ we can now deduce the corollaries of the main theorem.

\begin{corollary}{technical corollary}{} Let $(M,\m)$ be Haken with
complete and useful boundary pattern, and let $W$ be a union of
elements of $\uu{m}$. Then $\hscr(M,\m \rel W)$ is virtually
torsionfree and almost geometrically finite.
\marginwrite{technical corollary}
\end{corollary}

\begin{proof}{} By remark (i) above and proposition~\ref{relative
trick}, there is a complete and useful boundary pattern $\mtwo$ on $M$
so that $\hscr(M,\m\rel\abs{\mone})$ is isomorphic to a subgroup of
finite index in $\hscr(M,\mtwo)$. By theorem~\ref{main theorem},
$\hscr(M,\uu{m_2})$ is almost geometrically finite, hence so is
$\hscr(M,\uu{m}\rel\abs{\uu{m_1}})$.
\end{proof}

\begin{corollary}{rel W corollary}{} Let $M$ be a Haken 3-manifold and
let $W$ be a compact 2-dimensional submanifold of $\partial M$ such
that $\partial M-W$ is incompressible. Then $\hscr(M \rel W)$ is
virtually torsionfree and almost geometrically finite.
\marginwrite{rel W corollary}
\end{corollary}

\begin{proof}{} Let $\mone$ be the set of components of $W$, and let
$\m$ be the union of $\mone$ with the set of components of
$\ol{\partial M-W}$. The boundary pattern constucted above $\mtwo$
is complete, and by (ii) above it is useful. Since $\hscr(M\rel
W)\eq\hscr(M,\uu{m_2}\rel\abs{\uu{m_1}})$, corollary~\ref{technical
corollary} applies.
\end{proof}

\section [The Kontsevich Conjecture] {The Kontsevich Conjecture}
\label{kconj}\marginwrite{kconj}

Throughout this section, the symbol $f\cong g$ means that $f$ and $g$
are isotopic relative to the boundary of the manifold on which they
are defined. We first isolate a couple of technical steps that will be
needed in later arguments.

\begin{lemma}{homotopy implies isotopy}{} Let $M$ be a Haken manifold
containing an incompressible surface $G$, not necessarily connected.
Let $f$ and $g$ be two diffeomorphisms of $M$ which are homotopic
relative to $\partial M$. Then $f\cong g$. If $f$ and $g$ agree on
$G$, and the homotopy has trivial trace at a basepoint in each
component of $G$, then they are isotopic relative to $G\cup\partial
M$.
\marginwrite{homotopy implies isotopy}
\end{lemma}

\begin{proof}{} Replacing $f$ by $g^{-1}f$, we may assume that $f$ is
orientation-preserving, $g$ is the identity, and $f$ restricts to the
identity on $G$. As in theorem~II.6.1 of~\cite{Lau}, the homotopy can
be deformed relative to $M\times\partial I\cup \partial M\times I$ to
an isotopy. In particular, the traces at basepoints on components of
$G$ remain trivial. By theorem~\ref{Laudenbach}, we may assume that
the isotopy is relative to $G\cup\partial M$.
\end{proof}

\begin{lemma}{disconnected}{} Let $M$ be a compact 3-manifold, each
of whose components is Haken with incompressible boundary. Assume that
each component $M_0$ of $M$ has the property that if $g$ is a
diffeomorphism from $M_0$ to itself such that $g^n\cong 1_{M_0}$, then
there exists a diffeomorphism $h$ of $M_0$ such that
$g\cong h$ with $h^n\eq 1_{M_0}$. Then $M$ itself has this property.
\marginwrite{disconnected}
\end{lemma}

\begin{proof}{}
It suffices to consider the case when $g$ acts transitively on the
components of $M$. Let $M_1$ be a component and let $i$ be the
smallest positive integer such that $g^i$ preserves $M_1$. Since $M_1$
has the property, $g^i\cong h$ such that $h^{n/i}\eq 1_{M_1}$. Let
$M_2\eq g^{-1}(M_1)$. By isotopy relative to the boundary of $M_2$,
we may change $g\vert_{M_2}$ to $hg^{1-i}$, then $g$ has order $n$ on
each component of~$M$.
\end{proof}

We will need an extension of theorem~3 of \cite{H-T}.

\begin{proposition}{Heil-Tollefson torus boundary case}{} Let $M$ be
Haken with nonempty incompressible boundary. Assume that each
component of $\partial M$ is a torus. If $g$ is a diffeomorphism from
$M$ to itself such that $g^n\cong 1_M$, then $g\cong h$ such that
$h^n\eq 1_M$.
\marginwrite{Heil-Tollefson torus boundary case}
\end{proposition}

\begin{proof}{}
Consider the characteristic submanifold $\Sigma$ of
$(M,\ol{\uu{\emptyset}})$. Since $\partial M$ consists of tori, each
component of the frontier of $\Sigma$ is a torus, and is boundary
parallel exactly when it lies in a component of $\Sigma$ which is a
regular neighorhood of a component of $\partial M$. We induct on the
number of components of the frontier of $\Sigma$ that are not parallel
onto $\partial M$. If there are none, $M$ is either simple or a
Seifert fiber space. By theorem~3 of~\cite{H-T} and
lemma~\ref{homotopy implies isotopy}, the theorem is true for $M$. (In
the remainder of the proof, lemma~\ref{homotopy implies isotopy} must
be used in similar fashion to strengthen conclusions from
$\cite{H-T}$, but we will no longer mention these individually.)

We induct on the number of components of the frontier of $\Sigma$ that
are not parallel into $\partial M$. Since $\Sigma$ is unique up to
isotopy, we may assume that $g(\Sigma)\eq\Sigma$. Let $F$ be a
component of the frontier of $\Sigma$ such that $F$ is not parallel
into $\partial M$.  Let $\widehat{F}\eq\cup g^i(F)$, a collection of
components of the frontier of $\Sigma$. By induction and
lemma~\ref{disconnected}, it suffices to show that $g^n$ is isotopic
to $1_M$ relative to $F\cup\partial M$.

Fix an isotopy $H\colon g^n\cong 1_M$. The proof of lemma~9(ii)
of~\cite{H-T} contains most of the arguments needed to obtain our
conclusion, so we only explain the changes needed.  We refer to the
notation used there. The first paragraph of the proof (of Case~(ii))
in \cite{H-T} is not needed; since $\pi_1(M)$ is centerless the
condition $h(\gamma)\simeq\gamma$ holds automatically, as shown by the
argument for lemma~6 of~\cite{H-T}. For the next paragraph, we know
that $M$ cannot fiber over $S^1$ with fiber $F$, since $\partial M$ is
nonempty, therefore the first half of the paragraph shows that the
restriction of the homotopy to $\widehat{F}\times I$ is homotopic
relative to $\widehat{F}\times\partial I$ to a map into
$\widehat{F}\times I$. Therefore the homotopy at the start of the
third paragraph may also be assume to be relative to $\partial M$.
The remainder of the proof makes some rather delicate adjustments to
achieve that the restriction of $g$ to $\widehat{F}$ is periodic and
that the homotopy has trivial trace at a basepoint in each component
of $\widehat{F}$. Since these changes take place only in a regular
neighborhood of $\widehat{F}$ , the resulting homotopy is still
relative to $\partial M$. Then lemma~\ref{homotopy implies isotopy}
yields from this an isotopy relative to $F\cup\partial M$, to complete
the inductive step and the proof.
\end{proof}

\begin{theorem}{torsionfree mcg}{} Let $M$ be a Haken 3-manifold such
that $\partial M$ is nonempty and incompressible. Then $\hscr(M\rel
\partial M)$ is torsionfree.
\marginwrite{torsionfree mcg}
\end{theorem}

\begin{proof}{} Let $\langle g\rangle\in\hscr(M\rel \partial M)$ and
suppose that for $n>1$, $g^n\cong 1_M$. We must show that $g\cong
1_M$.

Let $T$ be the union of the torus boundary components of $M$. Suppose
first that $T\eq \partial M$. By proposition~\ref{Heil-Tollefson torus
boundary case}, $g\cong h$ with $h^n\eq 1_M$. Since $h$ is the
identity on $\partial M$, a theorem of M.\ H.\ A.\ Newman (see
proposition 3.1 of \cite{Lee}) shows that $h\eq 1_M$. Now suppose that
$T$ is not empty but $T\neq\partial M$. Form $N$ by gluing two copies
of $M$ together along $\partial M-T$, and let $D(g)$ be the
diffeomorphism of $N$ defined by taking $g$ on each copy of $M$. Since
$D(g)^n\cong 1_N$, the previous case shows that $D(g)\cong 1_N$. Let
$H\colon N\times I\to N$ be an isotopy from $D(g)$ to $1_N$. By
lemmas~7.2 and~7.3 of~\cite{Wald}, $H$ may be deformed to a homotopy
that preserves $G$. Therefore the trace of $H$ at a point in $G$ lies
in $G$. Since $D(g)$ is the identity on $G$, the trace is a central
element of $\pi_1(G)$. Since $G$ is not a torus, the center of
$\pi_1(G)$ is trivial. By theorem~\ref{Laudenbach}, $H$ may be
deformed to an isotopy relative to $G$. Repeating, we have an isotopy
relative to $\partial M-W$, so $g\cong 1_M$. This completes the case
when $T$ is not empty.

Now suppose that no component of $\partial M$ is a torus. Let $G$ be a
boundary component, and choose an essential simple closed curve
$\gamma$ in $G$. Let $G_1$ be a regular neighborhood of $\gamma$ in
$G$. Let $W$ be $S^1\times S^1\times I$, and let $G_2$ be a regular
neighborhood of $S^1\times \set{s_0}\times \set{0}$ in $S^1\times
S^1\times\set{0}$ for some $s_0\in S^1$. Form $N$ by identifying $G_1$
with $G_2$ and let $G_0$ be the incompressible surface in $N$ obtained
from $G_1$ and $G_2$. Since $G_0$ is incompressible in $M$ and $W$,
$N$ is Haken. Extend $g$ to a diffeomorphism $f$ of $N$ using the
identity on $W$. The isotopy $g^n\cong 1_M$ extends using the identity
on $W$ to an isotopy $f^n\cong 1_N$. Since $N$ has a torus boundary
component, the previous case implies that $f\cong 1_N$. By
theorem~\ref{Laudenbach}, $f\cong 1_N$ relative to $G_1$, and
therefore $g\cong 1_M$.
\end{proof}

Now we will weaken the hypothesis.

\begin{theorem}{torsionfree mcg rel F}{} Let $M$ be a Haken 3-manifold
and let $F$ be a nonempty compact 2-manifold in $\partial M$. Then
$\hscr(M\rel F)$ is torsionfree.
\marginwrite{torsionfree mcg rel F}
\end{theorem}

\begin{proof}{} Fix $\langle g\rangle\in\hscr(M\rel F)$ with $g^n$
isotopic to $1_M$ relative to $F$ for some $n>1$. We must prove that
$g$ is isotopic to $1_M$ relative to $F$.

Let $W$ be the union of the boundary components of $M$ that meet $F$.
On each component $X$ of $\ol{W-F}$, we have $g^n\vert_X\cong 1_X$.
Since $\partial X$ is nonempty, lemma~1.2 of~\cite{Hatcher-McCullough}
shows that $\hscr(X\rel\partial X)$ is torsionfree. Therefore
$g\vert_X\cong 1_X$, so we may change $g$ by isotopy relative to $F$
so that $g\vert_W\eq 1_W$. Also, $\pi_1(\Diff(X\rel\partial X))$ is
trivial, so we may assume that $g^n$ is isotopic to $1_M$ relative
to~$W$.

Form $N$ by attaching two copies of $M$ along $\partial M-W$. Then
$D(g)^n\cong 1_N$, so by theorem~\ref{torsionfree mcg}, $D(g)\cong
1_N$. As in the proof of theorem~\ref{torsionfree mcg}, the trace of
an isotopy from $D(g)$ to $1_N$ relative to $\partial N$ at each
component $G$ of $\partial M-W$ lies in $G$, so by
theorem~\ref{Laudenbach} we may deform the isotopy so that it
preserves $\partial M-W$. Therefore $g$ is isotopic to $1_M$ relative
to $W$ and hence relative to~$F$.
\end{proof}

Using corollary~\ref{rel W corollary}, we have the following immediate
consequence.

\begin{theorem}{Kontsevich Conjecture}{} Let $M$ be a Haken 3-manifold
with incompressible boundary, and let $F$ be a nonempty compact
2-manifold in $\partial M$, such that $\partial M-F$ is
incompressible. Then $\hscr(M\rel F)$ is geometrically finite.
\marginwrite{Kontsevich Conjecture}
\end{theorem}

\noindent From this we will obtain the following generalized version
of the Kontsevich Conjecture for Haken manifolds.

\begin{theorem}{Kontsevich Conjecture classifying spaces}{} Let $M$ be
a Haken 3-manifold with incompressible boundary, and let $F$ be a
nonempty compact 2-manifold in $\partial M$ such that $\partial M-F$
is incompressible. Then $\BDiff(M\rel F)$ has the homotopy type of a
finite complex.
\marginwrite{Kontsevich Conjecture classifying spaces}
\end{theorem}

\begin{proof}{} We will show that $\pi_i(\Diff(M\rel F))\eq 0$ for
$i\geq 1$. Since $\BDiff(M\rel F)$ is connected and
$\pi_{i+1}(\BDiff(M\rel F))\cong \pi_{i}(\Diff(M\rel F))$ for $i\geq 1$,
this implies that $\BDiff(M\rel F)$ is a $K(\hscr(M\rel
F),1)$-complex, so theorem~\ref{Kontsevich Conjecture} shows that
$\BDiff(M\rel F)$ has the homotopy type of a finite complex.

By the main theorem of \cite{Hat}, the homotopy groups
$\pi_i(\Diff(M\rel \partial M))$ vanish for $i\geq 1$, which gives the
assertion when $F\eq \partial M$. Otherwise, let $W\eq \ol{\partial
M-F}$. Restricting diffeomorphisms to $W$ is the projection map of a
fibration
$$\Diff(M\rel \partial M)\to\Diff(M\rel F)\to \Diff(W\rel \partial W)\ .$$

\noindent From the homotopy exact sequence of this fibration, we have
using \cite{Hat} again, that $\pi_i(\Diff(M\rel F))\cong
\pi_i(\Diff(W\rel \partial W))\eq 0$ for $i\geq 2$. We also obtain an
exact sequence
$$0\to \pi_1(\Diff(M\rel F))\to\pi_1(\Diff(W\rel\partial
W))\to\hscr(M\rel \partial M)\ .$$

\noindent
No component of $W$ is a 2-sphere, so elements in
$\pi_1(\Diff(W\rel\partial W))$ are classified by their traces
(nontrivial elements occur only for tori). Since the traces of an
isotopy from $1_M$ to $1_M$ at different basepoints are freely
homotopic, and $F$ is nonempty, all traces of an element of
$\pi_1(\Diff(M\rel F))$ must be  trivial in $\pi_1(M)$. Since $W$ is
incompressible, the restriction of an element of $\pi_1(\Diff(M\rel
F))$ to $W$ has trivial trace in each component of $W$, so is trivial
in $\pi_1(\Diff(W\rel\partial W))$. Therefore $\pi_1(\Diff(M\rel F))$
is trivial.
\end{proof}


\begin{thebibliography}{99}
{\footnotesize

\bibitem{B-S} A. Borel and J.-P.~Serre, Corners and arithmetic groups,
{\em Comment. Math. Helv.} 48 (1973) 436-491.

\bibitem{C-M} R. Canary and D. McCullough, Homotopy equivalences of
3-manifolds and deformation theory of Kleinian groups, in preparation.

\bibitem{C-V} L. S. Charlap and A. T. Vasquez, Compact flat Riemannian
manifolds III; the group of affinities, {\it Amer. J. Math} 95 (1973)
471-494.


\bibitem{Farkas-Kra} H. Farkas and I. Kra, {\it Riemann Surfaces,}
Springer-Verlag Graduate Texts in Mathematics~71 (1980).

\bibitem{F-L-P} A. Fathi, F. Laudenbach, and V. Poenaru, Travaux de
Thurston sur les surfaces, {\it Ast\'erisque} 66-67 (1979).

\bibitem{Gramain} A. Gramain, Le type d'homotopie du groupe des
diff\'eomorphismes d'une surface compacte, {\it Ann. scient. \'Ec.
Norm. Sup.} (4) 6 (1973) 53-66.

\bibitem{H-W} W. Hantzsche and W. Wendt, Drei dimensionale Euklidische
Raumformen, {\it Math. Ann.} 110 (1934) 593-611.

\bibitem{Harer1} J. Harer, The virtual cohomological dimension of the
mapping class group of an orientable surface, {\em Invent. Math.} 84
(1986) 157--176.

\bibitem{Harer2} J. Harer, The cohomology of the moduli space of
curves, in {\em Theory of Moduli,\/} ed. E. Sernesi, Springer-Verlag
Lecture Notes in Mathematics Vol. 1337 (1988) 138-221.

\bibitem{Harvey1} W. Harvey, Geometric structure of surface mapping
class groups, in {\it Homological Group Theory,} ed.~C. T. C. Wall,
London Math.~Soc.~Lecture Note Series No.~36 (1979), 255-269.

\bibitem{Harvey2} W. Harvey, Boundary structure of the modular group,
in {\it Riemann Surfaces and Related Topics: Proceedings of the 1978
Stony Brook Conference,} ed. I. Kra and B. Maskit, Annals of
Math.~Study No.~97 (1981), 245-251.

\bibitem{Hat} A. Hatcher, Homeomorphisms of sufficiently large
$P^2$-irreducible 3-manifolds, {\it Topology} 15 (1976) 343-347.

\bibitem{HatcherHarer} A. Hatcher, On triangulations of surfaces, {\em
Topology Appl.} 40 (1991) 189-194.

\bibitem{Hatcher-McCullough} A. Hatcher and D. McCullough, Finiteness
of classifying spaces of relative diffeomorphism groups of Haken
3-manifolds, preprint.

\bibitem{H-T} W. Heil and J. Tollefson, On Nielsen's theorem for
3-manifolds, {\it Yokohama Math. J.} 35 (1987) 1-20.

\bibitem{JS} W. Jaco and P. Shalen, Seifert fibered  spaces in
3-manifolds, {\it Mem. Amer. Math. Soc.} 220 (1979).

\bibitem{Joh} K. Johannson, {\it Homotopy Equivalences of 3-manifolds
with Boundary}, Springer-Verlag Lecture Notes in Mathematics 761
(1979).

\bibitem{KLR} Y. Kamishima, K. B. Lee and F. Raymond, The Seifert
construction and its application to infranilmanifolds, {\it Quart. J.
Math.} 34 (1983) 433-452.

\bibitem{K-P-S} A. Karrass, A. Pietrowski, and D. Solitar, Finite and
infinite cyclic extensions of free groups, {\it J. Australian
Math.~Soc.}~16 (1972), 458-466.

\bibitem{K} R. Kirby, Problems in low-dimensional topology, in {\it
Geometric Topology Part II} ed. W. Kazez, AMS/IP Studies in Advanced
Mathematics Vol. 2.2 (1997) 35-473.

\bibitem{Lee} J. S. Lee, Almost periodic homeomorphisms and p-adic
transformation groups on compact 3-manifolds, {\it Proc. Amer. Math. Soc.}
121 (1994) 267-273.

\bibitem{Lau} F. Laudenbach, Topologie de dimension trois. Homotopie et
isotopie, {\it Ast\'risque} 12 (1974) 11-152.

\bibitem{Lickorish} R. Lickorish, Homeomorphisms of non-orientable
two-manifolds, {\it Proc. Cambridge Phil. Soc.}~59 (1963), 307-317.

\bibitem{McCarty} G. S. McCarty, Homeotopy groups, {\it Trans. Amer.
Math. Soc.} 106 (1963) 293-304.

\bibitem{McCPoland} D. McCullough, Mappings of reducible 3-manifolds
in {\it Proceedings of the Se\-mes\-ter on Topology of the Stefan
Banach International Mathematical Center,} ed.~H.~Toru\'nczyk, Banach
Center Publications, Warsaw (1986), 61--76.


\bibitem{jdg} D. McCullough, Virtually geometrically finite mapping
class groups of 3-manifolds, {\it J. Diff. Geom.} 33 (1991) 1-65.

\bibitem{McCKorea} D. McCullough, {\it 3-Manifolds And Their
Mappings,} Global Analysis Center Lecture Note Series 26, Seoul
National University, 1995.

\bibitem{M-M} D. McCullough and A. Miller, Symmetric automorphisms of
free products, {\it Mem. Amer. Math. Soc.} 582 (1996), 1-97.

\bibitem{Orlik} P. Orlik, {\em Seifert Manifolds,} Springer-Verlag
Lecture Notes in Mathematics, Vol. 291 (1972).

\bibitem{Schneebeli} H. Schneebeli, On virtual properties and group
extensions, {\it Math. Z.} 159 (1978), 159-167.

\bibitem{Scottgeom} P. Scott, The geometries of 3-manifolds, {\it
Bull. London Math. Soc.} 15 (1983) 401-487.

\bibitem{VanMill} J. van Mill, {\em Infinite-Dimensional Topology,
Prerequisites and Introduction,} North-Holland Mathematical Library
Vol. 43 (1989).

\bibitem{W1} F. Waldhausen, Eine Klasse von 3-dimensionalen
Mannigfaltigkeiten I, II, {\em Invent. Math.} 3 (1967) 308--333,
{\bf 4}(1967) 87--117.

\bibitem{Wald} F. Waldhausen, On irreducible $3$-manifolds which are
sufficiently large, {\em Ann. of Math.} 87 (1968) 56--88.

\bibitem{West} J. West, Mapping Hilbert manifolds to ANR's: a solution
of a conjecture of Borsuk, {\em Ann. of Math.} (2) 106 (1977)J 1--18.

\bibitem{Zimmermann} B. Zimmermann, On the Hantzsche-Wendt manifold,
{\it Monatsh. Math.} 110 (1990) 321-327.

}
\end{thebibliography}
\end{document}